\documentclass[12pt]{amsart}

\usepackage[left=2.3cm, right=2.3cm, top=2.8cm, bottom=3.5cm]{geometry}
\usepackage{todonotes,ulem}
\def \Id{ {\mathrm {Id}}\,}
\usepackage{amssymb}
\usepackage{color,tikz,caption, pgfplots, amsmath,enumitem}
\usetikzlibrary{decorations.markings,arrows, calc}

\usepackage{wrapfig}
\def\bea#1\eea{\begin{align}#1\end{align}}
\usepackage{soul}
\setlist[enumerate]{itemsep=0mm}
\usepackage{framed}

\usepackage{amssymb,cancel}
\usepackage{mathrsfs}
\def \scr {\mathscr}
\definecolor{darkblue}{rgb}{0, 0, 0.5}
\usepackage{graphicx}
\usepackage[colorlinks=true, pdfstartview=FitV, linkcolor=darkblue, citecolor=blue, urlcolor=blue]{hyperref}
\definecolor{shadecolor}{rgb}{0.95, 0.95, 0.86}
\def \Sp{\mathfrak {Sp}}

\makeatletter
\@addtoreset{equation}{section}
\makeatother
\def\QED{\hfill $\blacksquare$\par\vskip 15pt}
\def\nn{\nonumber}
\newtheorem{theorem}{Theorem}[section]
\newtheorem{problem}[theorem]{RHP}

\newtheorem{coroll}[theorem]{Corollary}

\newtheorem{lemma}[theorem]{Lemma}

\newtheorem{proposition}[theorem]{Proposition} 
\newtheorem{definition}[theorem]{Definition}

\theoremstyle{remark}
\newtheorem{remark}[theorem]{Remark}

\def\br{\begin{remark}}
\def\er{\end{remark}}
\def\bt{\begin{theorem}}
\def\et{\end{theorem}}
\def\bc{\begin{coroll}}
\def\ec{\end{coroll}}
\def\bd{\begin{definition}}
\def\ed{\end{definition}}
\def\bp{\begin{proposition}}
\def \pa{\partial}
\def\ep{\end{proposition}}

\def\d{{\rm d}}

\def\wt{\widetilde}
\def\wh{\widehat}
\def\res{\mathop{\mathrm{res}}\limits_}
\numberwithin{equation}{section}

\def\le{\left}
\def\ri{\right}

\def\ov{\overline}
\def\1{{\bf 1}}

\def\bl{\begin{lemma}}
\def\el{\end{lemma}}
\def\be{\begin{equation}}
\def\ee{\end{equation}}

\def\C{{\mathbb C}}
\def \A{\mathbf A}
\def\l{\lambda}
\def\s{\sigma}
\def \B{\mathbf B}
\def\R{{\mathbb R}}
\def\N{{\mathbb N}}
\def\Z{{\mathbb Z}}
\def\G{\Gamma} 
\def\D{\Delta} 

\def\a{\alpha}

\def\e{\varepsilon}

\def\l{\lambda}

\def\r{\rho}
\def\s{ {\sigma}}

\def\x{\xi}
\def\z{\zeta}

\newcommand{\Iin}{J}
\newcommand{\Iex}{E}

\def\ra{\rightarrow}
\def\hf{\frac 12}
\newcommand{\CF}{\mathcal{F}}

\date{}

\begin{document}
%


\begin{flushright}
\end{flushright}
\vspace{0.2cm}
\begin{center}
\begin{Large}
On the spectral properties of the Hilbert transform operator on multi-intervals
\end{Large}

\bigskip
 M. Bertola$^{\dagger\ddagger\clubsuit}$\footnote{Work supported in part by the Natural
   Sciences and Engineering Research Council of Canada (NSERC)}\footnote{Marco.Bertola@\{sissa.it, concordia.ca\}},  
\\
A. Katsevich$^{\diamondsuit}$\\
A. Tovbis$^{\diamondsuit}$ \footnote{The work of AK and AT was partially supported by NSF grant DMS-1615124. The work of AK was supported also by NSF grant DMS-1906361. \{Alexander.Katsevich,Alexander.Tovbis\}@ucf.edu}
\\
\bigskip
\begin{small}
\begin{center}
\begin{enumerate}
\item [${\dagger}$] {\it  Department of Mathematics and
Statistics, Concordia University\\ 1455 de Maisonneuve W., Montr\'eal, Qu\'ebec,
Canada H3G 1M8} 
\item[${\ddagger}$] {\it SISSA/ISAS, via Bonomea 265, Trieste, Italy }
\item[${\clubsuit}$] {\it Centre de recherches math\'ematiques,
Universit\'e de Montr\'eal\\ C.~P.~6128, succ. centre ville, Montr\'eal,
Qu\'ebec, Canada H3C 3J7} 
\item [${\diamondsuit}$] {\it  University of Central Florida
	Department of Mathematics\\
	4000 Central Florida Blvd.
	P.O. Box 161364
	Orlando, FL 32816-1364
}\end{enumerate}
\end{center}
\end{small}
\end{center}
\bigskip
\begin{center}{\bf Abstract}\end{center}
Let $J,E\subset\R$ be two multi-intervals with non-intersecting interiors. Consider the following operator 
$$
A:\, L^2( J )\to L^2(E),\ (Af)(x) = \frac 1\pi\int_{ J } \frac {f(y)\d y}{x-y},
$$
and let $A^\dagger$ be its adjoint. We introduce a self-adjoint operator $\scr K$ acting on $L^2(E)\oplus L^2(J)$, whose off-diagonal blocks consist of $A$ and $A^\dagger$. In this paper we study the spectral properties of $\scr K$ and the operators $A^\dagger A$ and $A A^\dagger$. Our main tool is to obtain the resolvent of $\scr K$, which is denoted by $\scr R$, using an appropriate Riemann-Hilbert problem, and then compute the jump and poles of $\scr R$ in the spectral parameter $\l$. We show that the spectrum of $\scr K$ has an absolutely continuous component $[0,1]$ if and only if $J$ and $E$ have common endpoints, and its multiplicity equals to their number. If there are no common endpoints, the spectrum of $\scr K$ consists only of eigenvalues and $0$. If there are common endpoints, then $\scr K$ may have eigenvalues imbedded in the continuous spectrum, each of them has a finite multiplicity, and the eigenvalues may accumulate only at $0$. In all cases, $\scr K$ does not have a singular continuous spectrum. The spectral properties of $A^\dagger A$ and $A A^\dagger$, which are very similar to those of $\scr K$, are obtained as well.
\vspace{2mm}


\tableofcontents
%
%
\section{Introduction}


Let $J$, $E$ be two Lebesgue measurable subsets of $\R$. Consider the following operator 
\begin{align}
\label{Aoper}
A:\, L^2( J )\to L^2(E),\ (Af)(x) = \frac 1\pi\int_{ J } \frac {f(y)\d y}{x-y}, 
\end{align}
whose adjoint is 
\be\label{Adagger}
(A^\dagger g)(w) = \frac 1\pi \int_{E} \frac {g(x)\d x}{x-w}:\ L^2(E)\to L^2( J ).
\ee
An important and classical problem is to determine the nature of the spectrum of $A$, e.g., find its discrete and/or continuous parts and their multiplicities. When $\Iin=\Iex$, $A$ acts on the Hilbert space $L^2( J )$, and one can talk about the spectrum of $A$. In this setting the spectrum of $A$ for different sets $J$ was thoroughly studied starting in the 50's and 60's, see, e.g.,  \cite{koppinc59, kopp60, widom60, kopp64, pincus64, putnam65, rosenblum66}. For example, in the case where $J = E = \R$ the operator $-A$ is the usual Hilbert transform. The latter is well known to be anti-self-adjoint, and its spectrum consists of two eigenvalues $\pm i$. In particular, the spectrum of $A^\dagger A$ and $AA^\dagger$ is $+1$ (because the two operators are equal to the identity operator). This is easily seen by conjugating $A$ with the Fourier transform, which maps $A$ to the multiplication operator by $i {\rm sgn}(\xi)$, where $\xi$ is the Fourier variable and ${\rm sgn}$ is the signum function, see, for example, \cite{Kingv1}. Here and throughout the paper, the Fourier transform and its inverse are defined as follows:
\be\label{ft-inv}
\tilde\phi(\xi):=(\CF \phi)(\xi)=\frac1{\sqrt{2\pi}}\int_{ {\R}} \phi(t)e^{i\xi t}dt,\
\phi(t)=(\CF^{-1} \tilde\phi)(t)=\frac1{\sqrt{2\pi}}\int_{ {\R}} \tilde\phi(\xi)e^{-i\xi t}d\xi.
\ee
The operator $A$ is thus rather simple from the spectral point of view. In another known case \cite{OkE}, where $J = E$ is a finite interval, the operator $A$ in $L^2$ is not even a Fredholm operator (the range is dense, but not closed). In this case the spectrum of $A$ is absolutely continuous, of multiplicity 1, and coincides with the interval $[-i,i]$ of the imaginary axis \cite{koppinc59}.

More recently the problem was investigated in a number of new settings when $J\not=E$. Here, $J$ and $E$ can be intervals or multi-intervals (i.e., unions of finitely many non-intersecting closed intervals). When $J\not=E$, the spectral problem consists in the analysis of the operators $A^\dagger A$ and $AA^\dagger$. Such problems arise, for example, when solving the problem of image reconstruction from incomplete tomographic data, e.g. when solving the  interior problem of tomography \cite{yyww-07, yyw-07b, yyw-08, kcnd, cndk-08, aak14}. Different arrangements of $\Iin$ and $\Iex$ are possible, and they lead to different spectral properties of the associated operators. Spectral asymptotics for various arrangements of $\Iin$, $\Iex$, where each consists of a single interval  (the intervals can be disjoint or have a partial overlap), was obtained in \cite{kt12,  aak14, adk15}. In each of these cases the spectrum of the two operators is discrete. If $J$ and $E$ are disjoint, 0 is the only spectral accumulation point. If $J$ and $E$ overlap, there are two accumulation points: 0 and 1. When $J$ and $E$ are bounded intervals that touch at a point, the spectral set is $[0,1]$, and the spectrum is purely absolutely continuous with multiplicity one \cite{kt16}. Endpoints that are shared by both $\Iin$ and $\Iex$ are called {\it double} points. 
The analysis in \cite{kt12,  aak14, adk15, kt16} is based on the existence of  a differential operator that commutes with the finite Hilbert transform, which was found in \cite{kat10c, kat_11}.

Starting with \cite{BKT16}, the authors initiated the program of investigating the cases where $J$ and $E$ are multi-intervals subject to the restriction that the interiors of $J$ and $E$ are disjoint. 
In \cite{BKT16} we consider an arrangement, in which $E$ consists of two compact intervals, $J$ consists of any finite number of intervals that are all located between the two $E$ intervals, and  $\text{dist}(E,J)>0$. 
Since the use of commuting differential operator no longer applies when either $E$ or $J$ consists of more than one subinterval, the main tool in this paper is based on a matrix Riemann-Hilbert problem (RHP) approach to integral operators
with integrable kernels in the sense of \cite{IIKS}. 
The main findings of \cite{BKT16} include that the singular values of $A$ (ordered in decreasing order) tend to zero exponentially fast and an explicit expression for the leading term of the asymptotics. Let $\scr K$ be the self-adjoint operator acting on $L^2(E)\oplus L^2(J)$, whose off-diagonal blocks consist of $A$ and $A^\dagger$ (see \eqref{K-matr} below). In \cite{BKT16} we also showed that all the eigenvalues of $\scr K$ are simple and calculated the leading order behavior of its eigenvectors (in terms of Riemann Theta functions) as the spectral parameter $\l\to 0$. This operator is very convenient to work with as  $\scr K^2$
is a block diagonal operator with the blocks $AA^\dagger$ and $A^\dagger A$.
Speaking more generally, if $J$ and $E$ are arbitrary multi-intervals and $\text{dist}(E,J)>0$, it is easy to see that the operator $\scr K$ is compact, and the spectra of $A^\dagger A$ and $AA^\dagger$ are purely discrete. In fact, in the present paper we establish that $\scr K$  is of trace-class (that it is of Hilbert--Schmidt class is a simple exercise). 

The approach of  \cite{BKT16} works well when $\text{dist}(E,J)>0$, that is, when the integral operator $\scr K$ is
not singular. However, the case when $\text{dist}(E,J)=0$ leads to some technical difficulties, like, for example, construction of parametrices for the asymptotic solution of the RHP. These type of problems were overcome in 
\cite{BBKT}, where we use the RHP approach in the case where $J=[a,0]$ and $E=[0,b]$ for $a<0<b$, i.e. when $0$ is the only double point. 
The results of \cite{BBKT} match with and in some instances generalize those of  \cite{kt16}.
An arrangement where $J$ and $E$ have multiple common endpoints is considered in \cite{kbt19}. We assume there that $\Iin$ and $\Iex$ are multi-intervals, and their union is the whole line: $\Iin\cup\Iex=\R$. In this case, the corresponding RHP can be solved explicitly. Just as in \cite{BBKT}, the spectrum is the segment $[0,1]$, the spectrum is purely absolutely continuous, and its multiplicity equals to the number of double points. Additionally, in \cite{kbt19} we find an explicit diagonalization of the two operators.

In this paper we build on the results of \cite{BKT16} and \cite{kbt19} and extend the RHP approach further by allowing $J$ and $E$ to be general multi-intervals that can touch at any number of points, that is, $J$ and $E$ can have multiple double points. Our goal is to perform a qualitative analysis of the spectrum of $\scr K$ as well as of $A^\dagger A$ and $AA^\dagger$, which includes determining what spectral components it has and their multiplicities. It is quite interesting that without performing an explicit asymptotic analysis of the RHP when $\l\to 0$ that is similar to the one in \cite{BKT16}, and without access to an explicit solution of the RHP as in  \cite{kbt19}, much information can still be obtained by investigating $\scr K$ and the related RHP. Our main results are formulated in Section \ref{sec-main} below.
In addition to the RHP analysis, our main tools include 
the Kato-Rosenblum theorem on the stability of the absolutely continuous spectrum of a self-adjoint operator with respect to trace class perturbations.


\section{Main results} \label{sec-main}
As is stated in the introduction, the goal of the paper is to obtain the properties of the spectrum of $A^\dagger A$ and $AA^\dagger$ when  $E$ and $J$ are closed multi-intervals with disjoint interior $\displaystyle \mathop {J}^{\circ}   \cap \mathop{ E}^{\circ} = \varnothing$. Here and in what follows, $\displaystyle \mathop {U}^{\circ}$ denotes the interior of the set $U$.
The operator $A$ commutes with M\"obius transformations mapping $\R$ onto $\R$ (Lemma \ref{lemmamobius}) and, hence, it matters whether the sets $E,J$ have common points on the extended line. More precisely, if both $E$ and $J$ extend to infinity, we should consider $\infty$ as a common endpoint.  
The stated goal is essentially equivalent to studying the spectral properties of the self-adjoint operator $\mathscr K = A \oplus A^\dagger: L^2(U) \to L^2(U)$ with $U= E \cup J$. The latter is an operator with the kernel 
\be
\label{KKernel}
K(x,y) = \frac {\chi_{_J}(y)\chi_{_E}(x) -\chi_{_J}(x)\chi_{_E}(y)}{\pi(x-y)},
\ee
where, for a subset $U\subset \R$ we denote by $\chi_{_U}$ its indicator function. In matrix form, we can represent $\mathscr K$ as follows:
\be\label{K-matr}
\mathscr K=
\le[
\begin{array}{cc}
0 & A\\
A^\dagger & 0
\end{array}
\ri]: L^2(E)\oplus L^2(J) \to L^2(E)\oplus L^2(J).
\ee
The operator $\mathscr K$
is  a  convenient object to study because it is clearly self--adjoint and 
\be\label{K2-AK}
\scr K^2 = A A^\dagger \oplus A^\dagger A=
\le[
\begin{array}{cc}
A A^\dagger&0\\
0&A^\dagger A
\end{array}
\ri]:\,L^2(E)\oplus L^2(J) \to L^2(E)\oplus L^2(J). 
\ee
Thus, knowing $\Sp(\scr K)$, the spectrum of $\scr K$, it is easy to find the spectrum of $A A^\dagger$ and $A^\dagger A$. Hence, analysing $\Sp(\scr K)$ is also an important goal of this paper.

It is well known that
\be\label{SpK}
\Sp(\mathscr K)=\Sp_{ac}(\mathscr K)\cup\Sp_{sc}(\mathscr K)\cup\Sp_{p}(\mathscr K),
\ee
where $\Sp_{ac},\Sp_{sc}, \Sp_{p}$ denote the absolutely continuous, singular continuous, and  point spectra of $\mathscr K$, respectively. The main result of this paper, Theorem \ref{theo-main}, describes the connection between the geometry of 
the multi-intervals $E,J$ and the spectral components of $\mathscr K$. 
The main tool to studying the spectrum  of $\mathscr K$ is to construct the (nonsingular)  resolvent operator $\mathfrak R (\lambda) = \mathscr K (\Id-\frac 1 \lambda \mathscr K)^{-1} = (\Id-\frac 1 \lambda \mathscr K)^{-1}  -\Id  $.

It is a matter of inspection to ascertain that $\mathscr K$ is a Hilbert--Schmidt operator if $E$ and $J$ have no common endpoints (which implies that either $E$ or $J$ is bounded). In such case, therefore, the spectrum is purely discrete, and each eigenvalue has finite multiplicity. 

The term {\it endpoint} is used to denote an endpoint of any interval that makes up $E$ or $J$. An endpoint is called {\it simple} if it belongs only to one interval. An endpoint $z$ is called {\it double} if it belongs to two adjacent intervals of different types, that is, if 
$z\in E\cap J$. 
Naturally, two adjacent intervals of the {\it same} type are considered belonging to one interval. Our main theorem below
provides a detailed description of $\Sp(\mathscr K)$.

\bt
\label{theo-main} Let $\mathscr K = A \oplus A^\dagger: L^2(U) \to L^2(U)$ be the operator with the kernel \eqref{KKernel}. Here $U=E\cup J$, and $J,E\subset\R$ are multi-intervals with non-intersecting interiors. One has:
\begin{enumerate}
 \item $\Sp(\mathscr K)\subseteq [-1,1]$;
 \item There is an absolutely continuous component $\Sp_{ac}(\mathscr K) = [-1,1]$ of $\Sp(\mathscr K)$ if and only if there is a double
 point. Moreover, the multiplicity of $\Sp_{ac}(\mathscr K)$ is equal to the number of double points;
 \item  The end points $\l=\pm 1$ of the spectrum $[-1,1]$, as well as $\l=0$, are not eigenvalues. Moreover, $\mathscr K$ is of trace class if and only there are no double points. In this case, $\Sp(\mathscr K)$ consists only of eigenvalues and $\l=0$, which is the accumulation point of the eigenvalues;
 \item The eigenvalues of $\mathscr K$ are symmetric with respect to $\l=0$ and have finite multiplicities. Moreover, they can accumulate only at $\l=0$;
%
 %
 \item The singular continuous component is empty, i.e., $\Sp_{sc}(\mathscr K) = \varnothing$.
\end{enumerate}
\et
\br
According to assertion 1, all spectral components in \eqref{SpK} are subsets of $ [-1,1]$, i.e., the eigenvalues of $\mathscr K$ are embedded in the absolutely continuous spectrum  $\Sp_{ac}(\mathscr K)$ provided that both
 components are not empty.
 \er
\br When there is at least one double point, the presence of eigenvalues is not guaranteed. For example, it is shown in Proposition 4 of \cite{BBKT} that $A^\dagger A$ and $AA^\dagger$ do not have eigenvalues when $J=[b_L,0]$ and $E=[0,b_R]$ (i.e., $0$ is a double point). Here $b_L<0<b_R$. This implies that $\mathscr K$ does not have eigenvalues, because otherwise $\mathscr K^2$ (and $A^\dagger A$, $AA^\dagger$) would have had eigenvalues as well (see \eqref{K2-AK}).
\er
 
The proofs of assertions 1--3 are given in Section  \ref{sec-13}. They are based on the known facts about the spectrum of multi-interval Hilbert transforms,
see \cite{BKT16, kbt19, BBKT},
and the spectral trace class perturbation theorem of Kato-Rosenblum \cite{Kato}. An important part of our argument is Theorem \ref{ThmAtrace},
which states that the operator $A$ is of trace class provided that there are no double points.
The proofs of assertions 4 and 5 require a deep study of the solution $\G(z,\l)$ of a Riemann-Hilbert problem (RHP), which is associated with $\mathscr K$ on a certain  (infinite-sheeted) Riemann surface $\l\in\mathscr R$, see Section~\ref{sec-4}. We also add more details about $\G(z,\l)$ there.

Results from Theorem \ref{theo-main} can be naturally extended from the operator $\mathscr K$ to $\mathscr K^2$ (cf. \eqref{K2-AK}), thereby allowing us to obtain the analogue of Theorem~\ref{theo-main} for $A A^\dagger$ and $A^\dagger A$. 

\bt
\label{theo-mainK^2} Consider the operators $A$ and $A^\dagger$ defined by \eqref{Aoper} and \eqref{Adagger}, respectively. The operators $A A^\dagger$ and $A^\dagger A$ are unitarily equivalent. Also, the following assertions hold when $B=A A^\dagger$ or $A^\dagger A$:
\begin{enumerate}
 \item $\Sp(B)\subseteq [0,1]$;
 \item There is an absolutely continuous component $\Sp_{ac}(B) = [0,1]$ of $\Sp(B)$ if and only if there is a double point. The multiplicity of $\Sp_{ac}(B)$ is equal to the number of double points;
 \item $\l=0$, $\l= 1$ are not eigenvalues of $B$. Moreover, $B$ is of trace class if and only if there are no double points. In this case, $\Sp(B)$ consists only of eigenvalues and $\l=0$, which
 is the accumulation point of the eigenvalues;
\item The eigenvalues of $B$ have finite multiplicities, and they can accumulate only at $\l=0$;
%
 %
 \item The singular continuous spectrum of $B$ is empty, i.e., $\Sp_{sc}(B)=\varnothing$.
\end{enumerate}
\et

\section{Proof of Theorem \ref{theo-main}, assertions 1--3}\label{sec-13}
\subsection{Assertion 1: spectral interval.}
 Denote by $\mathcal H$ the Hilbert transform on $\R$. Then 
 $A=- \Pi_E \circ \mathcal H \circ \Pi_J$
  where $\Pi_J, \Pi_E$ are the projectors on $L^2( J ), L^2(E)$, respectively. It is well known that $\|\mathcal H\|=1$ (see e.g. \cite{Kingv1}, sec. 4.6), hence $\|A\|\leq 1$.  Consequently, the spectral radius of  $\mathscr K$ is also bounded by $1$, and assertion 1 is proven.
 
 \subsection{Assertion 2: absolutely continuous spectrum.}
 We have already commented that if there are no double points, i.e., the sets $ J , E$ are separated (in the extended line), then the self-adjoint operator $\mathscr K$ is a Hilbert--Schmidt operator, see also \cite{BKT16}. Then its spectrum is purely discrete and the eigenvalues (counted with multiplicity) form an $\ell^2$ sequence. In fact, it will be shown  below   that the operator $\mathscr K$ is of trace class if there are no  double points.
 
\subsection{Multiplicity of the continuous spectrum.}

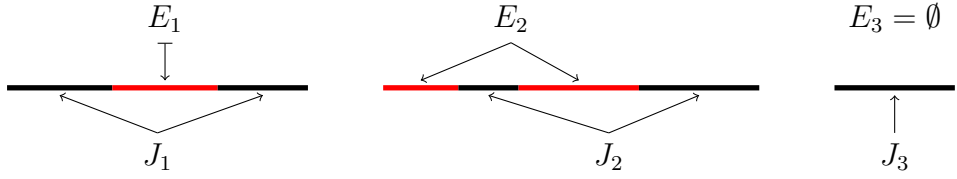
\begin{figure}
\begin{center}
\begin{tikzpicture}[scale=2]

\draw [line width =2] (-3,0) to (-2.3,0);
\draw [line width =2, red] (-2.3,0) to (-1.6,0);
\draw [line width =2] (-1.6,0) to (-1,0);
\node [above] at (-1.95, 0.3) {$ E_1 $};
\node [below] at (-2, -0.3) {$ J_1 $};
\draw [|->] (-1.95,0.3) to (-1.95, 0.05);
\draw [->] (-2,-0.3) to (-2.65,- 0.05);
\draw [->] (-2,-0.3) to (-1.3, -0.05);

\draw [line width =2,red ] (-0.5,0) to (-0,0);
\draw [line width =2] (0,0) to (0.4,0);
\draw [line width =2, red] (0.4,0) to (1.2,0);
\draw [line width =2] (1.2,0) to (2,0);
\node [above] at (0.35,0.3) {$ E_2 $};
\node [below] at (1,-0.3) {$ J_2 $};
\draw [->] (1,-0.3) to (0.2,- 0.05);
\draw [->] (1,-0.3) to (1.6, -0.05);

\draw [->] (0.35,0.3) to (-0.25, 0.05);
\draw [->] (0.35,0.3) to (0.8,0.05);

\draw [line width=2] (2.5,0) to (3.3,0);
\node [below] at (2.9,-0.3) {$ J_3 $};
\node [above] at (2.9,0.3) {$ E_3 = \emptyset $};
\draw [->] (2.9,-0.3) to (2.9,- 0.05);

\end{tikzpicture} 
\end{center}
\caption{An example of arrangement of $J$, $E$.}
\label{figIE}
\end{figure}

In this subsection we need a more detailed description of the multi-intervals $J,E$.
Let 
\be\label{JE-defs}
E = \bigcup_{j=1}^{ r } E_j,\qquad J = \bigcup_{j=1}^{ r }  J_j,
\ee
be the representations of $E$ and $J$ as unions of $r<\infty$ multi-intervals. Since the interiors of $E$ and $J$ do not intersect, we can arrange for the following properties to hold (see Figure \ref{figIE} for an illustration):
\begin{enumerate}
\item $J_i< J_j$ for $i<j$, i.e. the sets are ordered (ditto for the $E$ collection); 
\item ${\rm dist}(J_j, J_k)>0$ for $j\neq k$ (ditto for $E$);
\item for all $j=1\dots r $ the set $U_j := J_j\cup E_j$ is a single interval;
\item for all $i<j$ we have $E_i<  J_j$ and $J_i< E_j$ and in particular the distance ${\rm dist} (E_i, J_j)>0$ for $i\neq j$; 
\item For every $j=1,\dots,r$  the intersection  $J_j \cap E_j$ consists of $n_j$  endpoints of the sub-intervals.
\end{enumerate}

\bl
\label{lemmamobius}
Let $m(x) = \frac {ax + b}{cx +d}$ with $a,b,c,d\in \R$ and $ad-bc=1$ be a M\"obius  tranformation mapping $\R$ to $\R$; let $\mathcal U: L^2(\R, \d x) \to L^2(\R,\d x)$ be the corresponding unitary tranformation defined by:
\be\label{mathcal-U}
\mathcal U(f)(x) = \frac {f(m(x))}{ (cx + d)}. 
\ee
Then the Hilbert transform $\mathcal H$ on $\R$ commutes with $\mathcal U$: $\mathcal H \circ \mathcal U = \mathcal U \circ \mathcal H$.
\el
{\bf Proof.}
Let $g=\mathcal H f$. Then, using $dm/dx=\frac{1}{(cx+d)^2}$ we obtain
\be
\pi\frac{g(m(y))}{cy+d}=\int_\R\frac{f(\z)d\z}{(cy+d)(\z-m(y))}=\int_\R\frac{f(m(x))\frac{dx}{cx+d}}{(cy+d)(cx+d)(m(x)-m(y))}
=\int_\R\frac{f(m(x))\frac{dx}{cx+d}}{x-y},
\ee
where $\z=m(x)$. Note also that \eqref{mathcal-U} preserves the $L^2$ norm of $f$.
\QED

Lemma \ref{lemmamobius} implies that the spectral properties of $\mathcal H$ and all its possible restrictions are invariant under M\"obius transformations. In particular our operator $A$ is $A = -\Pi_E \mathcal H \Pi_J$, where $\Pi_U$ is the projection operator on the multi-interval $U$.  In this case $\wt A \,\mathcal U = \mathcal U A$, where $\wt A =-\Pi_{m(E)} \mathcal H \Pi_{m(J)}$.

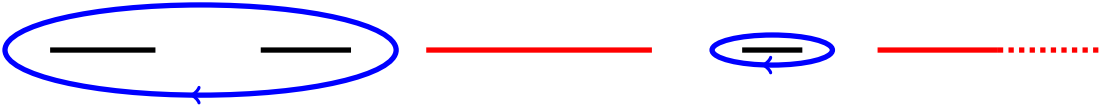
\begin{figure}
\begin{tikzpicture}[scale=2]

\draw [line width =2] (-3,0) to (-2.3,0);
\draw [line width =2] (-1.6,0) to (-1,0);

\draw [line width =2,red ] (-0.5,0) to (1,0);

\draw [line width =2] (1.6,0) to (2,0);

\draw [line width=2, red] (2.5,0) to (3.3,0);

\draw [line width=2, red, dotted] (3.3,0) to (4,0);

\draw [line width =2, blue,postaction={decorate,decoration={markings,mark=at position 0.75 with {\arrow[blue,line width=1.5pt]{<}}}}] (-2,0) ellipse (1.3 and 0.3);

\draw [line width =2, blue,postaction={decorate,decoration={markings,mark=at position 0.75 with {\arrow[blue,line width=1.5pt]{<}}}}] (1.8,0) ellipse (0.4 and 0.1);
\end{tikzpicture} 
\caption{An example of an arrangement of $J$, $E$ and $\gamma$ when there are no common endpoints. Note that $J$ or $E$ may have an unbounded component (in the picture it is $E$), but not both simultaneously. }
\label{FigIEgap}
\end{figure}

\bl
\label{ThmAtr0}
Let $A: L^2( J ) \to L^2(E)$  be the operator \eqref{Aoper} and $\mathscr K =A \oplus A^\dagger$. If $\text{dist}(J,E)>0$, then $\mathscr K$  is of trace-class.
\el
{\bf Proof.}
Let $\gamma$ be a union of Jordan curves separating $ J $ from $E$, see Figure \ref{FigIEgap}. Consider the Hilbert space $\mathcal L = L^2(U \cup \gamma, |\d z|)\simeq L^2( J )\oplus L^2(E) \oplus L^2(\gamma)$. 
Let $A_{\text{ext}}:\mathcal L\to\mathcal L$ be the operator with the kernel ${\chi_{_J}(y) \chi_{_E}(x)}/{(\pi (x-y))}$. Thus, $A_{\text{ext}}$ coincides with $A$ when the former is restricted to $L^2(J)\to L^2(E)$. We show that $A_{\text{ext}}$ is the product of two Hilbert--Schmidt operators, which immediately implies that $A$  is  of trace class. Indeed, let $T_k:\mathcal L\to\mathcal L,k=1,2$, be the following operators 
\be
T_1(f)(w) = \frac{\chi_\gamma(w)}{\pi}\int_J \frac {f(y)\d y}{w-y},\ 
T_2(g)(x) = \frac{\chi_E(x)}{2i\pi}\int_{\gamma} \frac {g(w)\d w}{w-x}.
\ee
The orientation of the contour $\gamma$ is chosen so that all points of $E$ are on the positive side. 
By construction, it follows immediately that both $T_{1}, T_2$ are Hilbert-Schmidt. Consider the composition
\be
T_2\circ T_1(f)(x) = \frac {\chi_{_E}(x)}{2i\pi^2}\int_{\gamma}\int_{J} \frac {f(y)\d y}{w-y} \frac {\d w}{w-x}.
\ee
An application of Cauchy's residue theorem shows that $T_2\circ T_1=A_{\text{ext}}$. Thus,  $A_{\text{ext}}$,  $A_{\text{ext}}^\dagger$, and $A_{\text{ext}}+A_{\text{ext}}^\dagger$ are all of trace class. Since $A_{\text{ext}}+A_{\text{ext}}^\dagger=\mathscr K\oplus\mathscr O$, where $\mathscr O: L^2(\gamma)\to L^2(\gamma)$ is the zero operator, we prove that $\mathscr K$ is of trace class.  \QED

\br
Iterating the argument in the proof of Lemma~\ref{ThmAtr0} one can represent $A_{\text{ext}}$ as a product of an arbitrary number of Hilbert--Schmidt operators. This means that the eigenvalues $\l_j$ (counted with multiplicity) of $\mathscr K$ form a sequence in $\ell^p$ for $\forall p \in(0,1]$, 
namely, 
\be
\sum_{j\geq 1} \lambda_j ^p<\infty,\ \ \ \forall p:\ \  0 <  p\leq 1. 
\ee
\er

\bl
\label{Lemma1} Suppose $r=1$ in \eqref{JE-defs}, i.e. $J=J_1$, $E=E_1$, and $U=U_1$.
Suppose $U = E\cup J$ is a single compact interval $[a,b]$, and $J$ and $E$ have $n$ endpoints in common. Then $\mathscr K:L^2(U)\to L^2(U)$ has absolutely continuous spectrum $[-1,1]$ of multiplicity $n$.
\el
{\bf Proof.}
There are two cases that need to be considered;
\begin{enumerate}
\item The leftmost and righmost sub-intervals in $U$ are of the same type: either both are parts of $J$ or both are parts of $E$;
\item The leftmost and rightmost sub-intervals in $U$ are of opposite types (e.g. the one on the left is a part of $J$, and the one on the right is a part of $E$).
\end{enumerate}
\paragraph{\bf First case.}
For definiteness suppose that both leftmost and rightmost intervals are part of $E$. 
Let $\mathscr K_{\text{{ext}}}: L^2(\R) \to L^2(\R)$ be the operator with the same kernel as $\mathscr K$ (see \eqref{KKernel}). The two operators act in a similar way, but $\mathscr K_{\text{ext}}$ acts on functions defined on all of $\R$.
Consider $\mathscr K_0: L^2(\R) \to L^2(\R)$ defined the same way as $\mathscr K_{\text{ext}}$ (i.e., with the kernel \eqref{KKernel}), but with $E$ replaced by $\wh E = E \cup U^c$. Here $U^c=\R\setminus U$, i.e. $\wh E$ ``extends'' $E$ to infinity.

The number $n$ of common endpoints between $J$ and $\wh E$ is the same as between $J$ and $E$. It is shown in Theorem~\ref{exact-soln}
that $\mathscr K_0$ has absolutely continuous spectrum $[-1,1]$ with multiplicity $n$. We also have 
\be\label{KOL}
\mathscr K_0 = \mathscr K_{\text{ext}} + \mathscr S,
\ee 
where $\mathscr S$ is the operator with the kernel 
\be
S(x,y) = \frac {\chi_{_{J}} (y) \chi_{_{U^c}} (x) -    \chi_{_{U^c}}(x) \chi_{_{J}}(y)}{\pi (x-y)}.
\ee
Since ${\rm dist}( J ,U^c)>0$, this operator is of trace class by Lemma~\ref{ThmAtr0} and, hence, $\mathscr K_0$ is a trace--class perturbation of $\mathscr K_{\text{ext}}$.  By the Kato-Rosenblum theorem \cite{Kato}, Theorem 4.4 of Chapter X, they have the same absolutely continuous spectrum with the same multiplicity. 

Finally, $\mathscr K_{\text{ext}}:L^2(\R)\to L^2(\R)$ coincides with the direct sum $\mathscr O \oplus \mathscr K$, where  $\mathscr K:\, L^2(U)\to L^2(U)$ is the original operator with the kernel \eqref{KKernel}, and $\mathscr O:L^2(U^c)\to L^2(U^c)$ is the zero operator. Therefore, the multiplicities of $\Sp_{ac}(\mathscr K)$ and $\Sp_{ac}(\mathscr K_{\text{ext}})$ are the same.

\paragraph{\bf Second case.}
Let $U =  J \cup E = [a,b]$. 
For definiteness suppose that the leftmost subinterval is part of $J$, and the other is part of $E$ so that $a\in  J , b\in E$. 
Define $\wh J$ and $\wh E$ by extending the corresponding two sub-intervals up to infinity. Let $\mathscr K_0$ be defined as $\mathscr K_{\text{ext}}$ with the replacements $ J \to \wh J $ and $E\to \wh E$. 
Let $E_\infty  := \overline{\wh E\setminus E} = [b,\infty)$  and $J_\infty=\overline{\wh J\setminus J} = (-\infty,a]$.
Similarly to the previous case  we have  
\be\label{KKS}
\mathscr K_0 = \mathscr K_{\text{ext}} +  \mathscr K_\infty +\mathscr S,
\ee
where the two operators $\mathscr S, \mathscr K_\infty$ have kernels, respectively,
\bea
S(x,y) &=  \frac {\chi_{_{J_\infty}} (y) \chi_{_E} (x) -    \chi_{_{E_\infty}}(x) \chi_{_J}(y)}{\pi (x-y)}
+  \frac {\chi_{_J} (y) \chi_{_{ E_\infty}} (x) -    \chi_{_E}(x) \chi_{_{J_\infty}}(y)}{\pi (x-y)}
\\
K_\infty(x,y) &= \frac {\chi_{_{J_\infty}} (y) \chi_{_{ E_\infty}} (x) -    \chi_{_{E_\infty}}(x) \chi_{_{J_\infty}}(y)}{\pi (x-y)}. 
\eea
Since ${\rm dist} (J_\infty, E) >0 $ and ${\rm dist}(E_\infty, J)>0$, it follows that $\mathscr S$ is of trace class as shown earlier.

However $\mathscr K_\infty$ is not trace-class because $E_\infty, J_\infty$ are both unbounded and ``meet'' at infinity. Indeed, by Lemma \ref{lemmamobius} the spectral properties of the Hilbert transform are invariant under M\"obius transformations that preserve the real line (i.e. $SL_2(\R)$).

Thus the spectral properties of $\mathscr K_\infty$ are equivalent to those of $\wh{\mathscr K}$ defined with $J = [-1,0]$ and $E=[0,1]$. This case was  analyzed in \cite{BBKT} where it was shown to have (only) absolutely continuous spectrum on $[-1,1]$ of multiplicity one.

On the other hand, Theorem \ref{exact-soln} from Appendix and Lemma \ref{lemmamobius}  show that $\Sp(\mathscr K_0)=[-1,1]$, with the absolutely continuous part of multiplicity $n +1$, where $n$ is the number of common endpoints between $J$ and $E$, and the 
additional  $+1$ multiplicity  is due to the fact that $\wh J$ and $\wh E$ meet at infinity (which can be mapped to a finite point by a M\"obius transformation). 

Since $\mathscr K_0$ is now  a trace--class perturbation of $\mathscr K_{\text{ext}}  + \mathscr K_\infty$ as per \eqref{KKS}, the multiplicity of $\Sp_{ac}(\mathscr K_0)$ must be the sum of the multiplicities of $\Sp_{ac}(\mathscr K)$ and $\Sp_{ac}(\mathscr K_\infty)$. The last statement follows, because the operator $\mathscr K_{\text{ext}}  + \mathscr K_\infty:L^2(\R)\to L^2(\R)$ coincides with the direct sum $\mathscr K \oplus \mathscr K_\infty$, where  $\mathscr K:\, L^2(E\cup J)\to L^2(E\cup J)$ and $\mathscr K_\infty:L^2(E_\infty\cup J_\infty)\to L^2(E_\infty\cup J_\infty)$ (we used here the same notation $\mathscr K_\infty$ for the original and restricted operators with a slight abuse of notation). It then follows that the multiplicity of $\Sp_{ac}(\mathscr K)$ equals $n$.
\QED
The following theorem completes the proof of assertion 2 and also the ``if and only if'' part of  assertion 3 of Theorem \ref{theo-main}.
\bt
\label{ThmAtrace}
Let $A: L^2( J ) \to L^2(E)$  be the operator \eqref{Aoper} and $\mathscr K =A \oplus A^\dagger$. 
Let  $n = \sum_{j=1}^{ r } n_j$ be the total number of  double endpoints in $U$, i.e., the total number of points  of contact between $J$ and $E$. If $n>0$, then the operator $\mathscr K$ has absolutely continuous spectrum $[-1,1]$ with multiplicity  $n$. 
\et
{\bf Proof.} 
Let $\mathscr K$ be the operator discussed above with the kernel \eqref{KKernel}. Consider the operators $\mathscr K_j:L^2(U_j)\to L^2(U_j)$ defined by the kernels 
\be\label{Kj-ker}
K_j(x,y)= \frac { \chi_{_{J_j}} (y) \chi_{_{E_j}} (x) -    \chi_{_{J_j}}(x) \chi_{_{E_j}}(y)} {\pi (x-y)}.
\ee
Consider also the operators $\mathscr H_{jk}: L^2(U) \to L^2(U)$ that are given by the kernels
\be
 H_{jk}(x,y) = \frac {\chi_{_{J_j}} (y) \chi_{_{E_k}} (x) -    \chi_{_{I_k}}(x) \chi_{_{E_j}}(y)}{\pi (x-y)},\ j\not=k.
\ee
Since ${\rm dist}(J_j, E_k)>0$ for $j\neq k$, all the operators $\mathscr H_{jk}$ are trace class by Lemma~\ref{ThmAtr0}. Using the two families of operators, represent the full operator $\mathscr K$ as follows
\be\label{main-K_split}
\mathscr K = \bigoplus _{j=1}^{ r } \mathscr K_j + \sum_{j<k} \mathscr H_{jk}.
\ee
Therefore $\mathscr K$ is a trace-class perturbation of the self-adjoint operator $\mathscr K_{\oplus}:=  \bigoplus _{j=1}^{ r }  \mathscr K_j$. The $\mathscr K_j$ are endomorphisms of the collection of  orthogonal subspaces  $\{L^2(U_j)\}_{j=1}^{ r }$ in $L^2(U)$.  The spectrum of $\mathscr K_{\oplus}$  is the disjoint union of the spectra of each $\mathscr K_j$. By Lemma \ref{Lemma1}, each $\mathscr K_j$ has absolutely continuous spectrum $[-1,1]$ with multiplicity $n_j$. Hence $\mathscr K_{\oplus}$ has absolutely continuous spectrum on $[-1,1]$ of multiplicity $n = \sum n_j$. By Theorem 4.4, p. 542 in \cite{Kato}, the absolutely continuous parts of $\mathscr K$ and $\mathscr K_{\oplus}$ are unitarily equivalent, and the theorem is proven.
\QED

 \subsection{\bf Assertion 3: point spectrum.}
 
We begin by proving that $\l=\pm 1$ are not eigenvalues of $\mathscr K$. Assume the opposite. Then, there exists $f\in L^2(U)$
  such that, for example,  $\mathscr K f=f$. Note that $|\mathscr K f(x)|=\,|\mathcal H f(x)|$ for $x\in U$, so  
  \be
  \|\mathscr K f\|_{L^2(U)} <   \|\mathcal H f\|_{L^2(\R)},
  \ee
  because $\mathcal H f$ is analytic in $\R\setminus U$.
  Since $ \|\mathcal H \|_{L^2(\R)}=1$, we obtain a contradiction
  \be
 \| f\|_{L^2(U)} = \|\mathscr K f\|_{L^2(U)} <   \|\mathcal H f\|_{L^2(\R)}\leq \|f\|_{L^2(U)}.
  \ee
  
Next, if $f\in L^2(J)$, then $Af\equiv 0$ if and only if $f\equiv 0$, since $Af$ is analytic in the  interior of $E$. 
Similarly, $A^\dagger g\not\equiv 0$ if $g\not\equiv 0$, where $g\in  L^2(E)$.
Therefore, $\l=0$ is not an eigenvalue of $\mathscr K$.
Finally, assume that there are no double points. According to Lemma \ref{ThmAtr0}, $\mathscr K$
is a trace class operator. But $\l=0$ is not its eigenvalue, therefore $\mathscr K$
must have a sequence of eigenvalues convergent to $\l=0$.
Thus, we proved assertion 3  of Theorem \ref{theo-main}.
%

\section{Proof of Theorem \ref{theo-main}, assertion 5}\label{sec-45}

To prove the remaining items of Theorem \ref{theo-main} we need to introduce the following RHP \ref{RHP-Gam} that is closely related
with the resolvent of $\mathscr K$. This approach goes back to \cite{IIKS}.
In the rest of the paper we will use the following  three Pauli matrices:
\be
\label{pauli}
\s_1 = \le[
\begin{array}{cc}
0&1\\
1&0
\end{array}
\ri]
\ \ \ \ 
\s_2 = \le[
\begin{array}{cc}
0&-i\\
i&0
\end{array}
\ri]
\ \ \ \ 
\s_3 = \le[
\begin{array}{cc}
1&0\\
0&-1
\end{array}
\ri].
\ee

\subsection{Riemann--Hilbert problem and the resolvent of $\mathscr K$ }\label{sec-4}
 \label{RHPRS}
 
Let us consider the following RHP.

\begin{problem}
 \label{RHP-Gam}
Find  a matrix--valued function $\Gamma(z;\lambda)$, such that for any
fixed   $\l\in\C\setminus{0}$, one has:
\begin{itemize}
\item[a)] the matrix $\Gamma(z;\lambda)$ is analytic together with its inverse 
in  $z\in\bar\C\setminus U$; 
\item[b)] $\Gamma(z;\lambda)$ satisfies the jump condition 
\be
\label{RHPGamma}
\Gamma_+ (z;\lambda) = \G_-(z;\lambda) \bigg(\1 - \frac {2i}\lambda f(z) g^T(z)\bigg),\quad  z\in  U=J \cup E, 
\ee
where 
\be\label{f,g}
f^T(x)= [ \chi_E(x),\chi_J(x)]; \ \ g^T(x) = [\chi_J(x), - \chi_E(x)];
\ee
\item[c)] $\G(\infty;\l) =\1$; and 
\item[d)] the limiting values  $\Gamma_\pm(z;\lambda)$ are in $L^2_{loc}$ near the endpoints of the intervals.
\end{itemize}
\end{problem}

The jump condition \eqref{RHPGamma} equivalently reads 
\bea
\G_+(z;\lambda)  &= \G_-(z;\lambda) \le[
\begin{array}{cc}
1 &0 \\
\frac {2i}\l & 1
\end{array}
\ri],\qquad z\in  J, 
\label{Ji}
\\
\G_+(z;\lambda)  &= \G_-(z;\lambda) \le[
\begin{array}{cc}
1 &-\frac {2i}\l \\
0 & 1
\end{array}
\ri],\qquad z\in E. 
\label{Je}
\eea

\br \label{rem-RHP-unique}
Using standard arguments, one can show that the requirement d) in the RHP \ref{RHP-Gam} implies the uniqueness of 
solution $\Gamma(z;\lambda)$ of the RHP \ref{RHP-Gam}. The existence of a solution in $\C\setminus [-1,1]$
will be proven in Theorem \ref{thm1} below.
Moreover, the solution has the following symmetries 
\be\label{symm}
\ov{\G(\ov z;\ov \l)} = \G(z;\l),\qquad 
\G(z;-\l) = \s_3 \G(z;\l)\s_3.
\ee
For example, the first symmetry follows 
 from the fact that  the matrix $V(z;\lambda) = \ov{\G(\ov z;\ov \l)} $ satisfies the same  RHP~\ref{RHP-Gam}. The proof of the second symmetry is also straightforward.
Additionally, it can be shown that if the solution to the RHP~\ref{RHP-Gam} exists, then it must satisfy
\be\label{det1}
\det \G(z;\l)\equiv 1.
\ee
\er

\br\label{rem-K}
One could guess that the RHP \ref{RHP-Gam} should be related with the operator $\mathscr K$, since
the kernel $K(x,y)$ of $\mathscr K$ given by \eqref{KKernel} can be represented as 
\be\label{Kfg}
K(x,y) = \frac {f^T(x) g(y)}{\pi(x- y)}.
\ee
\er

Let us now study the local  behavior of $\G(z;\l)$ near the endpoints.
Consider for example a simple right endpoint $z=a$ of $E$.   Denote by $\A$, $\B$ the first and  second  columns of the matrix $\G$, respectively. Then \eqref{Ji}-\eqref{Je} imply that  $\A(z)$ is analytic at $z=a$, and $\B(z) +\frac {2i}\l  \frac{\ln(z-a)}{2i\pi} \A(z)$ is analytic in the punctured neighborhood  of $z=a$.  The requirement d) of the RHP \ref{RHP-Gam} forces us to conclude that the latter expression is actually analytic (no pole). 
In other words, 
\be\label{simple-end}
\G(z;\l) = \mathcal O(1)\le[\begin{array}{cc}
1 &- \frac {2i}{\l} \frac{\ln(z-a)}{2i\pi}\\
0&1
\end{array}\ri].
\ee
Here and henceforth, $\mathcal O(1)$ denotes a matrix-valued function which is locally analytic in $z$ and invertible. A similar argument applies to all simple endpoints of $ J , E $.

Now consider a double endpoint. Without loss of generality we can place it at
 $z=0$ with $E $ locally on the  right of $z=0$, and $J$ - on the left. 
The first issue is the type of growth behavior that the entries of $\G$ have near $z=0$. To this end we observe that the jump matrices in \eqref{Ji}, \eqref{Je} are constant in $z$ and, therefore, we can analytically continue $\G(z)$ on the universal cover of a punctured neighborhood of $z=0$. Such analytic continuation has the following multivaluedness
\be
\G(z ) = \G(z{\rm e}^{2i\pi}) 
\le[
\begin{array}{cc}
1 & \frac {2i}\l\\
0 & 1
\end{array}
 \ri]
\le[
\begin{array}{cc}
1 & 0\\
\frac {2i}\l & 1
\end{array}
 \ri]
  =
 \G(z{\rm e}^{2i\pi}) 
\le[
\begin{array}{cc}
1 - \frac 4 {\l^2}  & \frac {2i}\l\\
\frac {2i}\l & 1 
\end{array}
 \ri] = \G(z{\rm e}^{2i\pi}) M_0,
 \label{GM0}
\ee
provided $\Im z>0$. Similar calculations show that $\G(z ) = \G(z{\rm e}^{2i\pi})M_0$ is valid for  $\Im z<0$ as well.

\begin{figure}
\begin{center}
\begin{tikzpicture}[scale=1.5]

\begin{scope}[xshift=-1cm]
\fill [white!90!black] circle(1.5);
\fill [white!80!black] (0,-1.5) arc (-90:90:1.5);
\draw[line width =2, red, postaction={decorate,decoration={markings,mark=at position 0.75 with {\arrow[black,line width=1.5pt]{>}}}}] (0,0.02) to (1,0.02); 
\draw[line width =2, blue, postaction={decorate,decoration={markings,mark=at position 0.25 with {\arrow[black,line width=1.5pt]{<}}}}] (0,-0.02) to (1,-0.02);

\node at (1,0) [below] {$\l=1$};
\node at (-1,0) [below] {$\l=-1$};
\draw[line width =2, green, postaction={decorate,decoration={markings,mark=at position 0.75 with {\arrow[black,line width=1.5pt]{>}}}}] (0,0.02) to (-1,0.02); 
\draw[line width =2,cyan, postaction={decorate,decoration={markings,mark=at position 0.25 with {\arrow[black,line width=1.5pt]{<}}}}] (0,-0.02) to (-1,-0.02); 
\end{scope}

\begin{scope}[xshift=2cm]

\node at (0.5,0) [right] {$\rho=\frac 1 2$};
\node at (-0.5,0) [left] {$\rho=-\frac 1 2$};
\fill [fill=white!90!black] (-0.5,-1) to (-0.5,1) to (0.5,1) to (0.5,-1);
\fill [fill=white!80!black] (-0.5,-1) to (-0.5,1) to (0,1) to (0,-1);

\draw[line width =2, red, postaction={decorate,decoration={markings,mark=at position 0.75 with {\arrow[black,line width=1.5pt]{<}}}}] (-0.5,0) to (-0.5,1); 
\draw[line width =2, blue, postaction={decorate,decoration={markings,mark=at position 0.75 with {\arrow[black,line width=1.5pt]{>}}}}] (-0.5,0) to (-0.5,-1); 

\draw[line width =2, green, postaction={decorate,decoration={markings,mark=at position 0.75 with {\arrow[black,line width=1.5pt]{<}}}}] (0.5,0) to (0.5,1); 
\draw[line width =2, cyan, postaction={decorate,decoration={markings,mark=at position 0.75 with {\arrow[black,line width=1.5pt]{>}}}}] (0.5,0) to (0.5,-1); 
\end{scope}
\end{tikzpicture} 
\end{center}
\caption{The slit $\l$-plane is mapped to the strip $|\Re \rho|<\frac 1 2$. The other strips $|\Re \rho-k|<\frac 1 2$ in the $\rho$-plane are mapped to the same slit $\l$-plane and represent the various sheets of the branched map $\l(\rho)$.}
\label{figrho}
\end{figure}
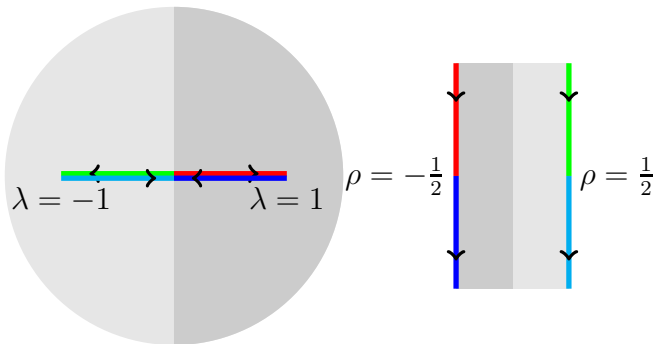

 Matrix $M_0$ 
 plays an important role in the analysis below. To calculate its eigenvalues and eigenvectors, it is convenient to
 introduce a new variable $\rho$, which is related to $\l$ as follows:
 \be
\label{rhosurf}
\rho(\lambda) =-\frac 1 2 +   \frac 1{ i\pi} \ln \le(\frac {1 - \sqrt{1-\lambda^2}}\lambda \ri), \qquad \l(\r)=-\frac{1}{\sin(\pi\r)},
\ee
We choose the branch of logarithm in \eqref{rhosurf} so that $\rho(\lambda) $ is a conformal mapping of $\C\setminus [-1,1]$
into the vertical strip  $|\Re (\rho)|<\frac 1 2$. We will also consider the analytic continuation  of this map 
as a map from the Riemann surface 
$\mathfrak  R$ of $\rho(\lambda)$ onto $\C$. 
 Figure \ref{figrho} provides a visualization
 of the map \eqref{rhosurf} 
 between the main sheet of $\mathfrak  R$, both shores of the branch cut $[-1,1]$ included, 
 and the vertical strip $|\Re\r|\leq \hf$.
 In general, each sheet of $\mathfrak  R$ is mapped onto the corresponding integer-shifted vertical strip $|\Re\r|\leq \hf$,
 so that $\r$ becomes a global coordinate on $\mathfrak  R$. 
Note that 
$\lambda(\r)$ is a single-valued meromorphic function on $\C$.  The
determination of the logarithm in \eqref{rhosurf} is consistent with
condition d) from RHP \ref{RHP-Gam}.

 Direct calculations show that $e^{\pm 2\pi i \r(\l)}$ are the eigenvalues of $M_0$, and  
\be\label{m-eig}
 C_{+} ^{-1}{\rm e}^{2i\pi \r\s_3} = M_0C_+^{-1},
\ee
where  
 \be\label{C_p} 
C_+(\r) :=  \left[ \begin {array}{cc}
1 & -e^{-i\pi\r}
\\
1  & e^{i\pi\r}
\end {array} \right] =  \left[ \begin {array}{cc}
-1 & \frac {i\l}{1 - \sqrt{1-\l^2}}
\\
-1  & \frac {i\l}{1+ \sqrt{1-\l^2}}
\end {array} \right], \qquad  \det C_+ = 2 \cos(\pi \rho).
\ee 
   To simplify notations, here and henceforth we often use $\r$ instead of $\r(\l)$.
  We also introduce
 \be\label{C-}
C_-(\rho) := C_+(\rho) \le[\begin {array}{cc}
1 & \frac {2i} \l\\
0 &1
\end{array}\ri]= \left[ \begin {array}{cc}
1 & -e^{i\pi\r}
\\
1  & e^{-i\pi\r}
\end {array} \right].
  \ee
 In the following proposition we derive the local behavior of $\G(z,\l)$ near a double point.


\bp
\label{propbeha}
Let $\l\in\C\setminus [-1,1]$. 
If $z=0$ is a double endpoint with $E $ adjacent to the right of $z=0$,
then any matrix valued function $\G(z;\l)$ satisfying  conditions a),b)  of the RHP  \ref{RHP-Gam} that is also $L^2_{loc}$  in a small disk $\mathbb D$
centered at $z=0$ can be written in the form 
\be\label{cond-double}
\G(z;\l) =Y(z;\l) z^{\rho(\l) \s_3} C_\pm (\r(\l)),\ z\in \mathbb D\cap \C^\pm,
\ee
where  $\C^\pm$ denote the upper/lower complex half plane and
$\rho(\lambda)$,    $C_\pm$ are given by \eqref{rhosurf}, \eqref{C_p}, and \eqref{C-} respectively.
Here $Y(z;\l)$ denotes a matrix valued function analytic near $z=0$ and with $\det Y(0;\l)\neq 0$.
\ep

\noindent  {\bf Proof.}
Let
\be
\label{parametrix}
P(z;\l):= z^{\rho  \s_3} C_\pm (\r).
\ee
We note that $\det P (z;\l) =2 \cos \pi\r$  is constant in $z$. 
A direct computation 
shows that $P$ satisfies the jump condition 
\bea\label{P-jumps}
P_+(z;\l) &= P_-(z;\l)\le[
\begin{array}{cc}
1 &\frac {-2i}{\l}\\
0&1
\end{array}\ri],\quad z\in \R^+,\\
P_+(z;\l) &= P_-(z;\l)\le[
\begin{array}{cc}
1 &0\\
\frac {2i}{\l}& 1
\end{array}\ri],\quad z\in \R^- \label{P-jumps1}.
\eea
Indeed, the first  equation follows from \eqref{C-}. The second  equation becomes ${\rm e}^{2i\pi \rho \s_3}C_+=C_+M_0$
as it  takes into account the jump of $ z^{\rho \s_3}$ on $\R^-$. 
Now \eqref{P-jumps1} follows from \eqref{m-eig}.

For $\l\not\in[-1,1]$, $\rho(\l)$ in \eqref{rhosurf} satisfies $|\Re \rho(\l)|<\frac 1 2 $, 
with $\Re \rho(\l)=\pm \hf$ being attained  on the $(0,1)$ and $(-1,0)$ parts of the branch cut  $(-1,1)$, respectively.
So, the required inequality is a consequence of the maximum principle for harmonic functions, see Figure \ref{figrho}. 
Hence  the matrix entries of $P(z;\l)$ are all in $L^2_{loc}$ near the origin for $\l\not \in [-1,1]$.

Now let $\Gamma(z;\r)$ satisfy  conditions a) and b)  of the RHP  \ref{RHP-Gam}
with entries in  $L^2_{loc}$ near $z=0$. Then $\G P^{-1}$ has no jumps in a neighborhood of the origin and, hence, it may only have an isolated singularity at $z=0$. The $L^2_{loc}$ condition together with $|\Re \rho(\l)|<\frac 1 2 $ implies that the singularity is removable. Thus the matrix $\Gamma$ has precisely the proposed representation \eqref{cond-double}. 
\QED

\br
\label{rem-Fuchs}
We should also point out that the solution  $\G(z;\l)$ of the RHP \ref{RHP-Gam}, if exists, 
solves a Fuchsian differential equation in $z$ of the form 
$$
\G'(z;\lambda) = \le(\sum_{z_j\in \pa J\cup \pa E} \frac {A_j}{z-z_j}\ri) \G(z;\lambda),
$$
 where the matrices $A_j$ are independent of $z$. These matrices, for a fixed $\lambda$, depend on the position of the endpoints according to the so--called Schlesinger equations \cite{JMU}, which express the fact that the monodromy representation induced by a fundamental solution of this ODE is independent of the position of the endpoints of the multi-intervals $J, E$. 
\er

The main tool for the analysis of the remaining assertions 4 and 5 of Theorem \ref{theo-main} is the following theorem for the so-called regularized resolvent defined by $\Id + \mathscr R(\lambda) = (\Id - \frac 1 \lambda \mathscr K)^{-1}$.
\bt
\label{thm1}
The resolvent $\mathscr R(\lambda)$ 
of $\mathscr K$ is an integral operator 
with the kernel 
\be
\label{resolGamma}
R(x,y;\lambda):=\frac 1 \lambda  \frac{f^T(x) \Gamma^{T}(x;\lambda) \Gamma^{-T}(y;\lambda) g(y)}{\pi (x-y)},\qquad 
\ee
where $\Gamma(z;\lambda)$ is the solution of the RHP \ref{RHP-Gam}.
%
Moreover, the operator $\Id-\frac 1\lambda \mathscr K$ has bounded inverse if and only if  the RHP 
 \ref{RHP-Gam} is solvable, and the solution is given by 
\be
\G(z;\l) = \1 - \int_{U} \frac{  F(x;\l) \cdot  g^T(x)\d x}{x-z}, 
\label{RestoRHP}
\ee
where $ F(x;\l) = (\Id  + \mathscr R(\l))[ f](x)$.
\et

{\bf Proof.} 
Assume that $ \G(z;\lambda)$ is the solution of the RHP \ref{RHP-Gam}. We first  show 
that the integral operator   $\mathscr R(\lambda)$ with the kernel \eqref{resolGamma} is bounded from  $L^2(U)$ into $L^2(U)$. Fix some $\l\in \C\setminus [-1,1]$.
We first note that, according to Remark \ref{rem-Fuchs},  $ \G_\pm(z;\lambda)$ is analytic
on $U$ with the exception of the  endpoints, where the local behavior of  $ \G_\pm(z;\lambda)$ is 
given either by \eqref{simple-end} (simple endpoint) or by Proposition \ref{propbeha} (double endpoint). Thus, the task of proving that $\mathscr  R$ is bounded requires only a local analysis in a neighborhood of each endpoint. If $z$ is a simple endpoint, the result is established in \cite{BKT16}. Even though the geometry of the intervals in \cite{BKT16} is slightly less general than the one here, the argument is purely local and applies in our situation as well.
Suppose now that $z=0$ is a double point. Since the problem is local, we can assume that $y$ is confined to a small neighborhood of $z=0$. Obviously, $\int_{U\setminus D_\e}R(x,y,\l)\phi(x)dx\in L^2(U)$,
where $D_\e$ is the $\e$- neighborhood of the origin, $\e>0$. Consider now the integral over $D_\e$.
Using the analyticity of $Y(z;\l)$ from Proposition \ref{propbeha}, we obtain
\be\label{G-inv-G}
 \Gamma^{-1}(y;\lambda) \Gamma(x;\lambda)=C^{-1}_+y^{-\r \s_3}\le[\1 +\mathcal O(x-y)\ri]x^{\r \s_3}C_+
\ee
uniformly in $x,y\in D_\e$ for a sufficiently small $\e$. Since the integral operator corresponding to 
the $\mathcal O(x-y)$ term is nonsingular and $|\Re \rho(\l)|<\frac 1 2$,
it remains to show that the integral operator
with the kernel
\be
R_0(x,y,\l)= \frac{ f^T(x)C^{-1}_+ \le(\frac xy \ri)^{\r \s_3}C_+  g(y)}{\l\pi (x-y)}
\ee
is a bounded operator in $L^2(U)$ (we have assumed that $0\in U$). 
According to \eqref{f,g},  the kernel $R_0$ is a linear combination of  
$\le(\frac xy \ri)^{\pm\r}$ and characteristic functions of $E,J$, so we can restrict our attention to
the integral operator 
\be
{\bf r}[\phi](y)= \int_U \frac{\le(\frac xy \ri)^{\r }\chi(x)\phi(x)}{x-y}dx,
\ee
where $\chi$ is either $\chi_{_E}$ or $\chi_{_J}$, and $\phi\in L^2(U)$. Using again that $|\Re \rho(\l)|<\frac 1 2 $ and appealing to Lemma 4.2 from \cite{GK}, p. 32, it is straightforward to conclude that ${\bf r}:L^2(U)\to L^2(U)$ is a bounded operator.
Thus, 
we proved that 
 $\mathscr R(\lambda)$ is a bounded operator in $L^2(U)$.

Let us now prove that the integral operator  $\mathscr R(\lambda)$ with the kernel \eqref{resolGamma} is the resolvent of   
$\mathscr   K$. 
The equation for the resolvent is 
\be
\label{resoll}
(\Id + \mathscr  
R)\circ\le(\Id - \frac 1 \l\mathscr   K\ri)  = \Id \ \ \Leftrightarrow\ \ \ 
 \mathscr  R -  \frac 1\l\mathscr  K = \frac 1 \l \mathscr  R \circ  \mathscr  K. 
\ee
As it was shown above, the kernel \eqref{resolGamma} defines a bounded integral  operator 
$\mathscr  R$ in $L^2(U)$, and we now verify that it satisfies \eqref{resoll}. Indeed,
the kernel $R \circ K$ of $ \mathscr  
R\circ\mathscr   K$ is
\bea
\frac 1 \l R \circ K(z,w) = \frac 1{( \l \pi)^{{2}}}  \int_{U} \frac { f^T(z) \Gamma^T(z;\l) \Gamma^{-T}(x;\l)   g(x)}{z -x}\frac { f^T(x)  g(w)}{x- w} \d x.
\label{RK}
\eea
Note that $\G^{-T}$ solves
\be
\G^{-T}_+ = \G^{-T}_-\le(\1 + \frac {2 i}\l  g f^T\ri),\quad z\in U,
\ee
so that for $x\in U$
\be
\Gamma^{-T}_{_+} - \Gamma^{-T}_{_-} = \frac{2i }\l \Gamma^{-T}_{_-} gf^T, 
\ee
and the right-hand side does not depend on the  side of the boundary  (recall that $ f^T  g \equiv 0$).
Thus \eqref{RK} yields
\bea\label{Prom1}
 &\int_{U}  f^T(z) \Gamma^T(z)\le(\Gamma^{-T}_+(x) - \Gamma^{-T}_-(x)\ri)   g(w) \frac {1}{(z-x)(x - w)} \frac{\d x}{2i\pi^{ 2}\l }
\cr
&=
 \int_{U}  f^T(z) \Gamma^T(z)\le(\Gamma^{-T}_+(x) - \Gamma^{-T}_-(x)\ri)   g(w) \frac 1{z - w} \le( 
 \frac 1{z - x}+  \frac 1{x - w}
 \ri) \frac{\d x}{2i\pi^{2} \l }.
\eea
To simplify notations, we drop the $\l$ dependence in $\G(x;\l)$ here and in the rest of the proof.
We show that \eqref{Prom1} splits into two essentially equal integrals. Indeed, choose $z\not\in U$  and, using Cauchy's theorem together with the fact that $\G(\infty)=\1$, we have  
\be\label{G^T}
\int_{U} \frac{\Gamma^{-T}_+(x) - \Gamma^{-T}_-(x)}{z - x} \frac{\d x}{2i\pi} = \1-\G^{-T}(z).
\ee
Substituting \eqref{G^T} into \eqref{Prom1} we finally obtain
\bea
&\frac 1{\l{\pi}(z - w)}  f^T(z) \Gamma^T(z)\Bigg[
  \overbrace{\int_{U}\frac {\le(\Gamma^{-T}_-(x) - \Gamma^{-T}_+(x)\ri)}{z - x}  \frac {\d x}{2i\pi}}^
  {\1-\Gamma^{-T} (z)}   
  +
 \overbrace{ \int_{U}\frac {\le(\Gamma^{-T}_-(x) - \Gamma^{-T}_+(x)\ri)}{ x - w}  \frac {\d x}{2i\pi}} ^{\Gamma^{-T}(w)-\1}
  \Bigg]
   g(w) \\
&= R(z,w) -\frac 1 \l K(z,w).
  \eea
  Thus we have shown that $\mathscr  R \mathscr  K = \mathscr  R -\frac 1 \l \mathscr  K$. Hence, the integral operator
   $\mathscr  R$ with the kernel $R$ given by \eqref{resolGamma} is the regularized resolvent.
  
 Vice versa, suppose now that $\Id -\frac 1 \l \mathscr  K$ is invertible. Define
 \be
  F(z) = \le(\Id -\frac 1 \l \mathscr  K\ri)^{-1} [ f],   
 \ee
 by which we mean the operator applied to each entry. 
 Then define 
 \be
 \label {418}
 \G(z;\l):= \1 - \int_{U}  \d x \frac {F(x) g^T(x)}{{\l\pi(x-z)}}.
 \ee
 We then  observe that $(\G_+-\G_-) f=0$ for  $x\in U$ so that we have 
 \bea
 \G_\pm(z;\l)  f(z) &=  f(z) -  \int_{U} \d x \frac { F(x)  g^T(x) f(z)}{{\l\pi(x-z)}} =  f(z) + 
 {\frac 1\l}\mathscr  K[ F] =  f + F - \le(\Id -{\frac 1\l}\mathscr  K\ri) [ F] \cr
 &=  f + F - \le(\Id -{\frac 1\l}\mathscr  K\ri) \bigg[\le(\Id -{\frac 1\l}\mathscr  K\ri)^{-1}[ f]\bigg]  =  f + F- f =  F.
 \eea
Thus,  for any $\l$ in the resolvent set of $\mathscr  K$ there exists $\G(z;\l)$ given by \eqref{418} 
that solves the RHP \ref{RHP-Gam}.
\QED

Using Remark \ref{rem-RHP-unique}, one can show that  the RHP \ref{RHP-Gam}
 is uniquely solvable if and only if the operator $\Id-\frac 1\lambda \mathscr K$ has a bounded inverse, i.e., when
 $\lambda\not =\Sp(\mathscr K)$, where $\Sp(\mathscr K)\subset[-1,1]$.

 \br\label{rem-R-pm}
 It is easy to show using the identity $f^T(z)g(z)\equiv 0$, $z\in U$, that the kernel $R(x,y;\l)$ in \eqref{resolGamma} does not have a jump across $U$. One can then combine this fact with the first equation of \eqref{symm} to prove that 
 $R(x,y;\bar\l)=\overline{R(x,y;\l)}$ when $x,y\in U$.
 \er 

\subsection{Study of the  spectrum of $\mathscr K$ by means of 
analytic continuation of the RHP solution across the spectral interval $(-1,1)$}

According to Theorems \ref{ThmAtrace}, \ref{thm1}, in the case of double points the RHP \ref{RHP-Gam} does  not  
have a solution for  any $\l\in[-1,1]$.  
In this subsection we discuss the meromorphic continuation of the solution   $\G(z;\l)$ to the RHP \ref{RHP-Gam}  over the segment 
$\l\in[-1,1]$ to the  Riemann surface $\mathfrak{R}$  of $\r(\l)$ beyond this segment. We will then use this continuation to analyze the 
resolvent $\mathscr R$ on $[-1,1]$.   Since $\r$ is a global coordinate on $\mathfrak R$, 
it will be convenient to introduce the notation
$\G(z;\r):=\G(z;\l(\r))$, where $\l(\r)$ is defined by \eqref{rhosurf}.

\begin{problem}
\label{genRHP}
Let $z_1,\dots, z_N$, $N\in\N$, be the  double endpoints (common endpoints $ J \cap E $). For a point 
$\r\in\C\setminus(\hf+\Z)$ we are 
looking for a matrix function  $\G(z;\r)$ with the following properties:
\begin{itemize}
\item[a)] the matrix $\Gamma(z;\r)$ is analytic together with its inverse 
in  $z\in\bar\C\setminus U$; 
\item[b)] $\Gamma(z;\r)$ satisfies the jump condition 
\bea
\G_+(z;\r)  = \G_-(z;\r) \le[
\begin{array}{cc}
1 &0 \\
 -2i\sin(\pi\r) & 1
\end{array}
\ri],\quad z\in  J, \\
\G_+(z;\r)  = \G_-(z;\r) \le[
\begin{array}{cc}
1 & 2i\sin(\pi\r) \\
0 & 1
\end{array}
\ri],\quad z\in E; 
\eea 
\item[c)] $\G(\infty;\r) =\1$; 
\item[d)] the limiting values  $\Gamma_\pm(z;\r)$ are in $L^2_{loc}$ for  any $z\in U\setminus\{z_1,\dots,z_N\}$; 
\item[e)] 
the local behavior of $\Gamma(z;\r)$
 near the double points is given by \eqref{cond-double} with $C_\pm$ given by \eqref{C_p}  and \eqref{C-}. 
%
%
\end{itemize}
\end{problem}

\bc\label{cor-RHPs}
For any $\r$ satisfying $|\Re\r|<\hf$ the solution $\G(z;\r)$ of the RHP \ref{genRHP} exists and  coincides 
with the solution $\G(z;\l)=\G(z;\l(\r))$ of the RHP \ref{RHP-Gam}.
\ec
{\bf Proof.}
According to \eqref{rhosurf}, 
conditions a), b) and  c)  of the RHPs \ref{RHP-Gam}  and \ref{genRHP} are the same. Moreover, conditions d) and e) of the 
 RHP \ref{genRHP} imply the condition d) of the  RHP \ref{RHP-Gam} provided $|\Re\r|<\hf$. Now the statement of the 
 corollary follows from Theorem \ref{thm1} and Proposition \ref{propbeha}.
\QED
%
We now aim to show that the solution $\G(z;\rho)$ of the RHP \ref{RHP-Gam} admits an extension to a meromorphic function of $\rho$ in  the whole $\rho$--plane. 
The proof proceeds in two steps:
\begin{itemize}
\item First, we prove that $\G$ admits an extension to a meromorphic function of $\rho$ in $\C\setminus (\Z+\frac 1 2)$. Observe that the points $\rho = \frac 1 2 +2\Z$ are all mapped to $\l = -1$, while the points $\rho = -\frac 12 +2\Z$ are mapped to $\l=1$. This implies that, in addition to $\r=\infty$, the poles    of $\G(z;\rho)$  can possibly accumulate at half integer $\r$;
\item we then prove
that near each of the points $\rho = \frac 1 2 +k$, $k\in\Z$,
$\G(z;\rho)$ is also meromorphic (i.e. has only finitely many poles).
\end{itemize}

The first point is proven in the next lemma.
\bl
\label{lemmaRS}
The solution of the RHP \ref{genRHP} admits a meromorphic extension to the domain  $\r\in\C \setminus (\Z+\frac 1 2) $.
\el

\noindent {\bf Proof.}
Let  $c$ be an  endpoint of $E$ or $J$, and let $\mathbb D_c$ be a small disk centered at $z=c$. We choose these disks centered at every endpoint of $E$ and $J$ small enough so that they are disjoint. 
Define $\Phi(z;\r) := \G(z;\r)$ outside of these disks, and
\begin{align}\label{Param}
\Phi(z;\r) &:= 
\le\{
\begin{array}{lc}
\G(z;\r)L(\pm(z-e);\r)^{\pm 1} 
& z\in \mathbb D_e
\\
\G(z;\r) U(\pm(z-f);\r)^{\pm 1}
 & z\in \mathbb D_f
 \\
 \G(z;\r) P_{_R}(z-q;\r)^{-1} & z\in \mathbb D_q\\
 \G(z;\r) P_{_L}(z-p;\r)^{-1} & z\in \mathbb D_p,\\
\end{array}
\ri.\\
\text{where}\quad
& L(z;\r):=  \le[\begin{array}{cc}1 & \frac{-\sin (\pi\r)\ln(z)}{\pi}\\ 0&1 \end{array}\ri],\ \ \ \ 
U(z;\r):= \le[\begin{array}{cc}1 & \\\frac{\sin (\pi\r)\ln(z)}{\pi}&1 \end{array}\ri],
\end{align}
$P_{_R}(z;\r)=P(z;\r)$ is the parametrix given by \eqref{parametrix}, and $P_{_L} = \s_2 P_{_R}\s_2$.  
Here $e$ is a simple  endpoint of $E$, $f$ is a simple endpoint of $J$,    $q$ is a double point having $E$ adjacent on the right, and $p$
is a double point having $E$ adjacent on the left, see Figure \ref{RHPphi}. {The sign `$+$' in \eqref{Param} is for the case that $e$ (respectively, $f$) is a right endpoint of $E$ (respectively, $J$), and the sign `$-$' -- for the left endpoints.

The results of Proposition \ref{propbeha} and the discussion immediately preceding it show that the matrix $\Phi(z;\r)$ is a piecewise analytic matrix-valued function on the complement of the contour $\Sigma$ that consists of the disks around the endpoints
together with the part of $U$ outside of these disks (see
 Figure \ref{RHPphi}), and $\Phi$ is uniformly bounded  with respect to  $z\in\C$. It is the solution of a RHP with 
 jumps on $\Sigma$, 
 where the jump matrices on the circles $\partial\mathbb D_c$ 
 depend {\it analytically} on $ \rho\in\C\setminus \Z+\frac 1 2 $. Moreover, the product of all the jump matrices at an intersection of any disk $\mathbb D_c$ with $U=E\cup J$ 
taken according to the orientation of each jump contour equals to the identity matrix $\1$.

Under these circumstances it is known, see Proposition 3.2 in \cite{FIKN} and also  the  original paper  \cite{Zhou1989},  that the solution $\Phi(z;\r)$ of the corresponding RHP either never exists 
or is meromorphic in $\rho$, with poles at  an exceptional locus of points that may accumulate only at the boundary of the domain of analyticity in the parameter space $\r$. 
The first option is not possible due to Corollary \ref{cor-RHPs}. Therefore, $\Phi(z;\r)$ is meromorphic with respect to  $\r\in\C\setminus (\Z+\frac 1 2)$, which implies the statement of the lemma.
\QED

We now need to prove that the solution $\G(z;\rho)$ of the RHP \ref{genRHP} is also meromorphic  
near $\rho = \frac 1 2 + k$, $k\in \Z$. To this end we first prove the lemma. 
\bl
\label{lemmaQk}
Let $z_j$, $j=1,\dots,N$, be a double endpoint with $E$ adjacent from the right. 
Define 
\be
\label{Qk}
Q_k(z,\rho):= \le\{
\begin{array}{cc}
\displaystyle
\le[
\begin{array}{cc}
1 & 0
\\
0 & \frac{{1}}{2\cos\pi\r}
\end{array}
\ri]
\le[
\begin{array}{cc}
1  & 0 
\\
-(z-z_j)^{2|k|-1} & 1
\end{array}
\ri] (z-z_j) ^{\rho \s_3} C_\pm (\r), 
& k= -1,-2,\dots
\\
\displaystyle
\le[
\begin{array}{cc}
 \frac{1}{2\cos\pi\r} & 0
\\
0 & 1
\end{array}
\ri]
\le[
\begin{array}{cc}
1  & -(z-z_j)^{2k+1}
\\
0 & 1
\end{array}
\ri] (z-z_j) ^{\rho \s_3} C_\pm (\r), & k=0,1,2\dots,
\end{array}
\ri.
\ee
where $\pm  \Im z>0$,  respectively.
Then the local behaviour of 
 the solution $\G(z;\r)$ near $z=z_j$ and 
 near the points $\r\in \hf+k$, where $k\in\Z$, can be represented by 
 \be\label{Gkj-right}
\G(z;\r) =Y_k^{(j)} (z;\rho)Q_k(z-z_j;\rho),\ j=1,\dots,N.
\ee  
Moreover, $\det Y_k^{(j)}\equiv {1}$,
and the matrix functions
 $Y_k^{(j)} (z;\rho)$ are analytic in $z$ in  $\r$-independent disks centered at  $z_j$.
If, instead of $E$, we have ${J}$ adjacent to $z_j$ from the right, then
the same statements hold provided we replace \eqref{Gkj-right} with
 \be\label{Gkj-left}
\G(z;\r) =Y_k^{(j)}(z;\rho)Q_k(z-z_j;\rho)\s_2.
\ee
\el
{\bf Proof.}
We consider the case when $E$ is adjacent to $z_j$ from the right.
The fact that $\det Y_k^{(j)}\equiv 1$ follows from \eqref{C_p}.
The analyticity of $Y_k^{(j)}(z;\rho)$ with respect to $z$ in a $\r$-independent neighborhood of the double point $z_j$ follows from the fact that $\G(z;\rho)Q_k^{-1}(z-z_j;\rho)$ is analytic near $z_j$.  The other case when $J$ is adjacent to $z_j$ from the right can be considered analogously.
\QED

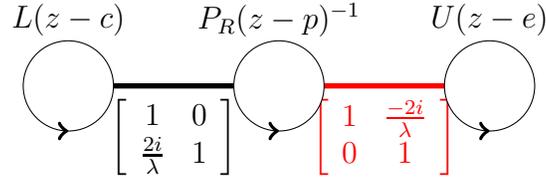
\begin{figure}
\begin{center}
\begin{tikzpicture}[scale=4]

\draw [line width =2] (-3,0) to node[below] {$\le[\begin{array}{cc}
1 &0 \\
\frac {2i}\l & 1
\end{array}
\ri]$} (-2.3,0);
\draw [line width =2, red] (-2.3,0) to  node [below]{$\le[\begin{array}{cc}
1 &\frac {-2i}\l \\
0 & 1
\end{array}
\ri]$}(-1.6,0);
\foreach \x\U in {-3/ ,-2.3/ ,-1.6/  }
{
\draw[ fill=white, postaction={decorate,decoration={markings,mark=at position 0.75 with {\arrow[black,line width=1.5pt]{>}}}}] (\x,0) circle[radius=0.15]; 
}

\node at (-3,0.22) {$L(z-c)$};

\node at (-2.3,0.22) {$P_R(z-p)^{-1}$};
\node at (-1.6,0.22) {$U(z-e)$};

\end{tikzpicture} 
\end{center}
\caption{An example of the  contour $\Sigma$ supporting the jump discontinuities of the RHP for the matrix $\Phi$. Indicated near each arc is the corresponding jump matrix. The black line segment represents $J$, and the red line segment represents $E$.}
\label{RHPphi}
\end{figure}

\bp\label{prop-cond-l=pm1}
Assume $z_j=0$ in \eqref{Qk}.
Then the matrices $Q_k(z;\rho)$ with $z\neq 0$ are analytic in some neighborhoods of  $\rho = \frac 1  2+k$ for all $k\in\Z$.
 For $\r=\hf+k$, where $k=-1,-2,\dots$, we have 
\be \label{beha+1}
\lim_{\e\ra 0}Q_k(z;\hf-|k|+\e)=\le[
\begin{array}{cc}
z^{\hf-|k|} &
(-1)^k iz^{\hf-|k|} \\
\frac {(-1)^{k}  z^{-\hf+|k|}\ln z}{\pi} 
& z^{-\hf+|k|} {-\frac {1}{\pi i} }  z^{-\hf+|k|} \ln z 
\end{array}
\ri],
\ee
where the expression is for $z$ in the upper half plane. Similarly  for $\rho=\frac 1 2 + k$ and $k=0,1,\dots$ we have 
\be
\lim_{\e\ra 0}Q_k(z;\hf+k+\e)=\le[
\begin{array}{cc}
\frac {(-1)^{k}  z^{\hf+k}\ln z}{\pi} 
& -z^{\hf+k} {+\frac {1}{\pi i} }  z^{\hf+k} \ln z 
\\
z^{-\hf-k} &
(-1)^k iz^{-\hf-k} \\
\end{array}
\ri].
\ee
These equations can be easily modified for $\Im z<0$ using \eqref{C-} and \eqref{Qk}.
\ep

{\bf Proof.}
Multiplication of all the factors in the first case of  \eqref{Qk} yields
\be \label{prod1}
Q_k(z;\r)=\le[
\begin{array}{cc}
z^\r &
-{\rm e}^{-i\pi\r}z^\r \\
\frac {{-}z^{\r+2|k|-1}{+} z^{-\r}}{{2}\cos\pi\r} 
&{-} \frac {z^{\r+2|k|-1}{\rm e}^{-i\pi\r} +z^{-\r}{\rm e}^{i\pi\r}}{2\cos\pi\r} 
\end{array}
\ri].
\ee
Substituting $\r=\hf+k+\e$ into \eqref{prod1} we derive \eqref{beha+1} after some algebra.  All other cases can be considered analogously.
\QED
With these preparations we can finally prove the meromorphic continuation of $\G(z;\r)$ onto the whole $\rho$--plane.
\bl
\label{propRS}
The solution of the RHP \ref{genRHP} admits an extension to a meromorphic function of $\r\in\C$.
\el
\noindent {\bf Proof.}
We know from Lemma \ref{lemmaRS} that the matrix $\G(z;\rho)$ admits a meromorphic extension to $\r\in \C \setminus (\Z+\frac 1 2)$. In principle, that lemma does not rule out  an accumulation of poles near the points $\Z+\frac 1 2$. 
Therefore we still need to prove that $\G(z;\rho)$ is meromorphic also in  neighborhoods of each of the points 
$\rho_k = \frac 12 + k$, $k\in \Z$.
This part of the  proof is now only a minor revision of Lemma \ref{lemmaRS} and, therefore, we use the same notation from  
that proof.  Fix $k\in \Z$ and  define 
\begin{align}\label{Paramk}
\Phi(z;\r) &= 
\le\{
\begin{array}{lc}
\G(z;\r)L(\pm (z-e))^{\pm 1} 
& z\in \mathbb D_e
\\
\G(z;\r) U(\pm (z-f))^{\pm 1}
 & z\in \mathbb D_f
 \\
 \G(z;\r) Q_{k}(z-q;\r)^{-1} & z\in \mathbb D_q\\
 \G(z;\r) \s_2  Q_{k}(z-p;\r)^{-1}\s_2 & z\in \mathbb D_p,\\
\end{array}
\ri.
\end{align}
where $L, U$ are defined in \eqref{Param}, and $Q_k$ -- in \eqref{Qk}. 
Define $\Phi(z;\r) = \G(z;\r)$ outside of the disks 
from \eqref{Paramk}.
 Here, as in  Lemma \ref{lemmaRS}, $e$ is a simple endpoint of $E$,  $f$ is  a simple  endpoint of $J$, $q$ is a double point having $E$ adjacent on the right, whereas $p$ is a double point having $E$ adjacent on the left.
The sign $+$ is for the case when  $e$ ($f$, respectively) is a right endpoint of $E$ ($J$ respectively) and the sign $-$ for the left endpoints.
 %
   {Choose a small neighborhood $S$ of $\rho = \hf +k$.
The same reasoning used in Lemma \ref{lemmaRS} now applies 
to $\r\in S$ due
to Lemma \ref{lemmaQk}. Thus,  we conclude that $\G(z;\r)$ is meromorphic in a neighborhood of $\rho = \hf +k$. 
\QED

 \subsection{\bf Absence of singular continuous  spectrum, assertion 5.}

Before proceeding we briefly summarize the consequences of  Lemma \ref{propRS}, see also  \eqref{rhosurf} and Figure \ref{figrho}.
Since  the main strip $\Re \rho\in (-\hf,\hf)$ corresponds to $\l\not\in [-1,1]$, Theorem \ref{thm1} implies  that none of the poles of $\G(z;\rho)$ (which is the solution to the RHP \ref{genRHP})  occurs in that strip.
Thus, if any, the only poles in the closure of the main strip  may occur on the lines $\Re \rho = \pm \frac 1 2$, that is,
on the shores of the segment $(-1,1)$ of the spectral $\l$ plane. If $\mathscr K$ has a continuous spectrum, 
that is, if there is at least one double point, then the poles of $\G(z;\rho)$, as we are going to show in Section \ref{sec-4},  
correspond to the embedded point spectrum   of  $\mathscr K$. Thus, to complete the spectral description of $\scr K$, in the following theorem we
prove the absence of singular continuous spectrum. This will prove assertion 5 of Theorem \ref{theo-main}.

\bt\label{theo-sing-spec}
The singular continuous component of $\Sp(\mathscr K)$ is empty, i.e., $\Sp_{sc}(\mathscr K) = \varnothing$.
\et

\noindent {\bf Proof.} Let $\G(z,\l)$ be the solution of the RHP \ref{RHP-Gam}. 
It is clear that a pole $\r_0,\, |\Re \r_0|=\hf$, of the solution $\G(z;\rho)$ of the RHP \ref{genRHP} corresponds 
to the pole $\l_0=\l(\r_0)$, $\l_0 \in(-1,1)$ of $\G(z,\l)$. 
Since $|\G(z,\l)|\equiv 1$, the poles (in $\l$) of 
$\G(z,\l)$ and $\G^{-1}(z,\l)$ coincide. Thus, the kernel $R(x,y,\l)$ of the resolvent operator of $\mathscr K$, given by 
\eqref{resolGamma}, is meromorphic in $\l$. Then so is the jump $\D_\l R(x,y,\l):= R(x,y,\l_+)-R(x,y,\l_-)$ over the 
interval $\l\in(-1,1)$.
In particular, $\D_\l R(x,y,\l)$ has no more than finitely many poles on any closed subinterval of $(-1,1)\setminus \{0\}$.

Pick any $f\in C_0^{\infty}(\mathring{J}\cup \mathring{E} )$. Let $f_1=\mathscr P f$ be the projection of $f$ onto the direct sum of all the eigenspaces of $\mathscr K$ (i.e., the subspace of discontinuity with respect to $\mathscr K$, see e.g. Section X.1.1 in \cite{Kato}). Set $f_2:=f-f_1$. Let $\mathscr E_\l$ denote the resolution of the identity associated with $\scr K$. Using the properties of $\mathscr E_\l$  (i.e., $\mathscr P^2=\mathscr P$ and $\mathscr E_\l\mathscr P=\mathscr P \mathscr E_\l$, see, e.g., Sections VI.5.1 and X.1.1 in \cite{Kato}), we have by starting with $\mathscr Pf_1=f_1$:
\be\label{sigma}
\sigma_{f_2}(\l):=(\mathscr E_\l(f-f_1),f-f_1)=(\mathscr E_\l f,f)-(\mathscr E_\l f_1,f_1). 
\ee
We want to prove that $\sigma_{f_2}(\l)$ is smooth for any $\l\not=0$ not in the {point} spectrum of $\mathscr K$.
Let $[\l_1,\l_2]$ be any interval that does not contain any eigenvalue of $\mathscr K$ such that $\l_1<\l<\l_2$. Without loss of generality we may assume $\l_1>0$. The case $\l_2<0$ can be considered analogously.
Then $\mathscr E_\l f_1=\mathscr E_{\l_1} f_1$, and the second term on the right in \eqref{sigma} is locally constant with respect to $\l$. 
Also, $\mathscr E_\l=\mathscr E_{\l_1}+\mathscr E_{[\l_1,\l]}$. According to \cite{DS57} p. 921, $\mathscr E_\l$ is computed by the formula 
\begin{equation}\label{ResOfIdGeneral}
\mathscr{E}_{\l}=\frac{-1}{2\pi i}\lim_{\epsilon\to0^+}\int_{-\infty}^{\l}\left[\mathscr{R}_1(t+i\epsilon)-\mathscr{R}_1(t-i\epsilon)\right]\d t,
\end{equation}
if $\l$ is not an eigenvalue. 
Here $\mathscr R_1(\l) := \le(\l\Id - \mathscr K\ri)^{-1}$. Clearly, $\mathscr R_1(\l)=(1/\l)(\Id+\mathscr R(\l))$.
Therefore,
\begin{equation}
(\mathscr E_{[\l_1,\l]} f,f)=\frac{-1}{2\pi i\l}\int_{\l_1}^{\l}\int_U\int_U\Delta_t R(x,y,t)f(x)\bar f(y)\d x\d y\d t
\end{equation}
is a locally smooth function of $\l$ because $R(x,y,t)$ is $C^\infty$ on $\text{supp}f\times \text{supp}f\times [\l_1,\l_2]$. By construction, $\sigma_{f_2}(\l)$ is a continuous function of $\l$. We just proved that it may fail to be smooth only at the eigenvalues of $\mathscr K$ and at $\l=0$. Since the number of eigenvalues of $\mathscr K$ is finite in any set $[-1,-\e)\cup(\e,1]$, $\epsilon>0$, this implies that $\sigma_{f_2}(\l)$ is absolutely continuous. Since the span of $C_0^{\infty}
(\mathring J\cup \mathring E )$
 is dense in $L^2(U)$, we see that $\mathscr K $ has no singular continuous spectrum. 
\QED

\section{Proof of Theorem \ref{theo-main}, assertion 4, and Theorem~\ref{theo-mainK^2}}\label{sec-4}

\subsection{Proof of Theorem \ref{theo-main}, assertion 4}

We will now prove the remaining part of assertion 4 from Theorem \ref{theo-main}, namely, 
that each eigenvalue $\l_0\in (-1,1)$ of $\mathscr K$ has a finite dimensional eigenspace. 
The symmetry of the eigenvalues with respect to $\l=0$ follows by noticing that if $\mathscr K(f,g)^T=\l(f,g)^T$, $\l\not=0$, $f\in L^2(E)$, $g\in L^2(J)$, $(f,g)\not\equiv0$, 
then $\mathscr K(-f,g)^T=-\l(-f,g)^T$, see \eqref{K-matr}.
This also follows from the symmetry \eqref{symm} of the solution $\G(z;\l)$ of the RHP \ref{RHP-Gam}.}

\bp\label{prop-pole} 
Any pole of the solution  $\G(z;\l)$ to the RHP \ref{RHP-Gam} is a simple pole.
\ep
\noindent {\bf Proof.}
Since $\mathscr K$ is a bounded, self-adjoint operator, the resolvent $\mathscr R$ of $\mathscr K$ has only simple poles, and  
$|\l|^{-1}\| \Id+\mathscr R(\l)\|=[\text{dist}(\l,\Sp(\mathscr K))]^{-1}$ (see \cite{ReedSimon1}, Example 2 on p. 224).
So, if $\l_0$ is a pole of $\mathscr R(\l)$, then $\| \mathscr R(\l)\|\leq \frac{c}{|\l-\l_0|}$, $\l\to\l_0$, 
for some $c>0$. Then, according to \eqref{RestoRHP}, $\G(z;\l)$ also has a pole at $\l_0$ whose order can not exceed one.
\QED

Let $\l_0\in (-1,1)$ be a simple pole of $\G(z;\l)$ with  the Laurent expansion near $\l_0$ given by 
\be\label{G-exp-pole}
\G(z;\l)= \frac { \G^{0}(z)}{\l-\l_0 } + \G^1(z) + \mathcal O(\l-\l_0),
\ee
where the term $\mathcal O(\l-\l_0)$ is uniform near any point of analyticity (in $z$) of $ \G(z;\l)$.
The representation \eqref{G-exp-pole} can be modified in a natural way so that it works near simple and double endpoints, see \eqref{simple-end} and Proposition \ref{propbeha}
respectively.
%
%

\bp\label{prop-const}
The matrix $\G^{0}(z)$ in \eqref{G-exp-pole} can be written as follows
\be
\label{G0}
\G^0(z) =  \le[\begin{array}{c} 
a\\ b
\end{array}\ri] \Psi(z),
\ee
where $a,b\in \C$ are constants that are not both zero, and
the vector $ \Psi(z) := [\psi_1(z), \psi_2(z)]$ has the jump condition and asymptotics given by
\be
\label{jumpG0_a}
\Psi_+(z)= 
\Psi_-(z) \le[\1 - \frac {2i}{\l_0} \s_+\chi_{_E}(z) + \frac {2i}{\l_0} \s_-\chi_{_J}(z)  \ri], \quad
\Psi(z) = \mathcal O(z^{-k}) \ \ \text{as}~~z\to\infty  
\ee
 with some $k=1,2,\dots$. Here
\be
\label{s_pm}
\s_+ := \le[
\begin{array}{cc}
0&1\\
0&0
\end{array}
\ri],\qquad 
\s_- := \le[
\begin{array}{cc}
0&0\\
1&0
\end{array}
\ri].
\ee
\ep
{\bf Proof.}
The jump conditions in \eqref{jumpG0_a} follow immediately from \eqref{G0}, \eqref{G-exp-pole} and the RHP \ref{RHP-Gam}. 
Also, note that the matrix $\G^0(z)$ is analytic at $z=\infty$ and vanishes because $\G(z;\l) = \1 + \mathcal O(z^{-1})$. This  implies that there is $k\in \N$ such that $\G^{0}(z) = \mathcal O(z^{-k})$, and this implies also  the same bound for $\Psi$ in \eqref{jumpG0_a}.

Thus, it remains to prove \eqref{G0}.
Using \eqref{G-exp-pole} and the relation
 $\G^{-1}(z;\l) = \s_2 \G^T(z;\l) \s_2$, we obtain
\begin{align}
&\frac 1{\l_0 (\l-\l_0)^2} \frac { g(y)^T \s_2 (\G^0)^T(y) \s_2 \G^0(x) f(x)}{(x-y)} + 
\nn\\ 
&+\frac 1{ (\l-\l_0)}\bigg[ 
 \frac { g(y)^T \le( \s_2 (\G^0)^T(y) \s_2 \G^1(x) + \s_2 (\G^1)^T(y) \s_2 \G^0(x)\ri) f(x)}{\l_0 (x-y)} 
 +\nn\\
 &\hspace{2.3cm}-\frac { g(y)^T \s_2 (\G^0)^T(y) \s_2 \G^0(x) f(x)}{\l_0^2 (x-y)} 
\bigg]+\label{1stpole}
\mathcal O(1).
\end{align}  
The numerator of the second order pole equals
\bea\label{G0calc}
\le[
\G^0_{22}(y)\chi_J(y) + \G^0_{21}(y) \chi_E(y)\ ,\ 
-\G^0_{12}(y)\chi_J(y) - \G^0_{11}(y) \chi_E(y)
\ri]
\le[\begin{array}{c}
\G^0_{11}(x) \chi_E(x) +\G^0_{12}(x) \chi_J(x)\\
\G^0_{21}(x) \chi_E(x) + \G^0_{22}(x) \chi_J(x)\\
\end{array}\ri]
\cr
=
 \bigg(
\G^0_{22}(y)\G^0_{11}(x)- \G^0_{12}(y)\G^0_{21}(x)
\bigg) \chi_J(y) \chi_E(x) 
+ \bigg(
\G^0_{21}(y)\G^0_{12}(x)- \G^0_{11}(y)\G^0_{22}(x)
\bigg) \chi_E(y) \chi_J(x) 
\cr
+ \bigg(
\G^0_{22}(y)\G^0_{12}(x)- \G^0_{12}(y)\G^0_{22}(x)
\bigg) \chi_J(y) \chi_J(x) 
+ \bigg(
\G^0_{21}(y)\G^0_{11}(x)- \G^0_{11}(y)\G^0_{21}(x)
\bigg) \chi_E(y) \chi_E(x). 
\eea
If follows from Proposition \ref{prop-pole} that the kernel $R$ should have a a first order pole in $\l$,
so that \eqref{G0calc} must be identically zero. 
This expression is identically zero if and only if the two rows of $\G^0$ are proportional by a {\it constant}. This is so because, for example, $\G^0_{22}(y)\G^0_{12}(x)- \G^0_{12}(y)\G^0_{22}(x)\equiv 0$ implies $\G^0_{22}(y)/ \G^0_{12}(y)=\G^0_{22}(x)/\G^0_{12}(x)$ and both sides must be constants because they depend on different variables.
\QED

\bc\label{cor-res}
The leading order term $P(x,y;\l)=\frac { P_{\l_0}(x,y)}{\l-\l_0}$ of $ R(x,y;\l )$ near a pole $\l=\l_0$ does  not have a jump across $U$, where 
\be\label{pole-R}
P_{\l_0}(x,y):=\frac { g(y)^T \le( \s_2 (\G^0)^T(y) \s_2 \G^1(x) + \s_2 (\G^1)^T(y) \s_2 \G^0(x)\ri) f(x)}{\l_0 (x-y)}.
\ee
\ec

{\bf Proof.}
Equation \eqref{pole-R} follows directly from \eqref{1stpole}, where the  $O((\l-\l_0)^{-2})$ term is zero,
and \eqref{G0}. Indeed, from \eqref{G0} it follows that $(\G^0)^T(y) \s_2 \G^0(x) \equiv 0$ and hence the last term in the simple--pole term of \eqref{1stpole} is zero.
Also, since $P(x,y;\l)$ is the leading order term of $R(x,y;\l)$, it follows that $P_{\l_0}(x,y)$ has no jump across $U$.
\QED


\bl\label{cor-fin-dim}\label{lem-G12}
The kernel $P_{\l_0}(x,y)$ is degenerate. 
\el

{\bf Proof.}
Substituting \eqref{G-exp-pole} into  the RHP \ref{RHP-Gam} we obtain the following jump conditions and the asymptotics
for the Laurent coefficients 
 $\G^{0,1}(z)$: 
\bea
\label{jumpG0}
\G^{0}(z)_+&= 
\G^{0}(z)_- \le[\1 - \frac {2i}{\l_0} \s_+\chi_{_E}(z) + \frac {2i}{\l_0} \s_-\chi_{_J}(z)  \ri],\\
\G^0 (z) &= \mathcal O(z^{-1}), \ \ z\to\infty,  \ \ \\ 
\label{RHPG1}
\G^1(z)_+ &= \G^1(z)_- \le[\1 - \frac {2i}{\l_0} \s_+\chi_{_E}(z) + \frac {2i}{\l_0} \s_-\chi_{_J}(z)  \ri] 
 + \G^{0}(z)_- \le[
  \frac {2i}{\l_0^2} \s_+\chi_{_E}(z) - \frac {2i}{\l_0^2} \s_-\chi_{_J}(z) 
 \ri], \\
 \G^{1}(z) &= \1  + \mathcal O(z^{-1}),\label{G-asympt}
\eea
where 
\be
\det \G^0\equiv 0 ,\qquad {\rm Tr}\bigg(\G^0 \s_2 (\G^1)^T\s_2\bigg)\equiv 0.
\label{tracecond}
\ee
The determinant and trace conditions  follow from the property $\det \G(z;\l)\equiv 1$ (see \eqref{det1}).

%
Inserting \eqref{G0} into  the trace condition \eqref{tracecond} gives
\be
0\equiv [b,-a]\G^1(z)  \le[\begin{array}{c}
\psi_2(z)\\-\psi_1(z)
\end{array}\ri] \ \ \Rightarrow \ \  \ 
[b,-a]\G^1(z) = h(z)  \Psi(z) \label{abeq}
\ee
for some scalar function $h(z)$ to be identified. Next, our goal is to show that $h(z)$ is a rational function. This is done in three steps.
\begin{itemize}
 \item
First, 
multiplying \eqref{RHPG1} on the left by $[b,-a]$ and noticing that $[b,-a]\G^0(z)\equiv 0$, we obtain that   $h_+(z) = h_-(z)$ for $z\in U$. 
This means that $h(z)$ extends analytically across $U$. 
\item
Second, let 
$z_0\in\C$ be any point other than an endpoint of $J$ or $E$  where {\it both} components of $ \Psi$ vanish. 
Then $\G^0(z_0)=0$ and so, according to \eqref{G-exp-pole}, $\det \G(z_0;\l)=\det \G^1(z_0) +\mathcal O(\l-\l_0)$.
If a zero  $z_0$   of the vector $ \Psi$ has multiplicity $\mu$, 
then $h(z)$ may have a pole of order at most $\mu$ since the left side of \eqref{abeq} is bounded.

\item
Zeroes of $\Psi$ cannot accumulate at any  $z_*\in \C$. Assuming the opposite, let $z^*$ be a point where the zeroes of $\Psi$ accumulate. If $z_*$ is not an endpoint, then $\G(z;\l)$ and, consequently, $\G^0(z)$ are analytic at $z_*$ (see \eqref{G-exp-pole}). By \eqref{G0}, such accumulation implies that $\G^0(z) \equiv  0$, and $\l_0$ is not a pole.
Suppose $z_*$ is an endpoint, for example, $z_*$ is a double point. Since $z^{\rho(\l) \s_3} C_\pm (\r(\l))$ is analytic in $\l$ in a neighborhood of $\l_0$, then \eqref{cond-double} implies that $Y(z;\l)$ has a pole at $\l_0$. Thus, we can repeat the same argument with the matrix function $Y(z;\l)$ (which is analytic near $z_*$)
 and its residue 
$Y^0(z)=\G^0(z)z^{-\r(\l_0)\s_3}C^{-1}_\pm(\r(\l_0))$. In the case when $z_*$ is a simple endpoint one can use \eqref{simple-end} instead of \eqref{cond-double}.
\item
Finally, observe that by \eqref{G-asympt} the left hand side of the second equation in \eqref{abeq} tends to $[b,-a]$ at $z=\infty$.
Hence, we conclude that $h(z)$ has polynomial growth of  degree not exceeding $ k$ (see \eqref{jumpG0_a}) and, therefore, 
 according to  Liouville's theorem,
 it is a rational function.
\end{itemize}

We also observe  that  \eqref{abeq} and \eqref{G0} imply
\be
\label{166}
\s_2(\G^0)^T(y) \s_2 \G^1(x) = h(x) \le[
\begin{array}{c}
\psi_2(y)\\
-\psi_1(y)
\end{array}
\ri][ \psi_1(x), \psi_2(x)].
\ee


Substituting \eqref{166} into \eqref{pole-R} we obtain
\bea\label{res-form}
P_{\l_0}(x,y)&= 
 \frac { g(y)^T \le( \s_2 (\G^0)^T(y) \s_2 \G^1(x) + \s_2 (\G^1)^T(y) \s_2 \G^0(x)\ri) f(x)}{\l_0 (x-y)} \cr
& = i\frac { h(x) -h(y)} {\l_0(x-y)} { g(y)^T \s_2  \Psi^T(y) \Psi(x) f(x)}.
\eea
The expression $\frac { h(x) -h(y)} {\l_0(x-y)}$ is a finite linear combination of products of rational functions in $x$ and $y$ separately with at most as many terms as the degree of the scalar rational function $h(z)$. Thus $P_{\l_0}$ is a degenerate kernel. 
\QED

To obtain a more explicit expression for $P_{\l_0}(x,y)$ we simplify \eqref{res-form}. Suppose $h(x)=S_t(x)/S_b(x)$, where $S_t$ and $S_b$ are some polynomials without common factors. Then
\be\label{bezout}
\frac { h(x) -h(y)}{x-y}=\frac {S_t(x)S_b(y) -S_t(y)S_b(x)}{(x-y)S_b(x)S_b(y)}
=\frac{\sum_{m,n} b_{mn} x^m y^n}{S_b(x)S_b(y)},\ b_{mn}=b_{nm},
\ee
where $(b_{mn})$ is the B\'ezout matrix of the polynomials $S_t(x),S_b(x)$. Using \eqref{pauli} and \eqref{f,g}, we get
\be\label{P_factor}\begin{split}
g(y)^T \s_2  \Psi^T(y) \Psi(x) f(x)
&=\begin{bmatrix}\chi_{_J}(y)& - \chi_{_E}(y)\end{bmatrix}
\le[
\begin{array}{cc}
0&-i\\
i&0
\end{array}
\ri]
\begin{bmatrix} \psi_1(y)\psi_1(x) & \psi_1(y)\psi_2(x)\\
\psi_2(y)\psi_1(x) & \psi_2(y)\psi_2(x)\end{bmatrix}
\begin{bmatrix} \chi_{_E}(x) \\ \chi_{_J}(x)\end{bmatrix}\\
&=-i(\psi_1(x)\chi_{_E}(x) +\psi_2(x) \chi_{_J}(x))(\psi_1(y)\chi_{_E}(y) +\psi_2(y) \chi_{_J}(y)).
\end{split}
\ee
Combining \eqref{res-form}--\eqref{P_factor} gives
\be\label{P-symm}\begin{split}
P_{\l_0}(x,y) =& B(x,y) H(x)H(y),\\
B(x,y):=&\l_0^{-1}\sum_{m,n} b_{mn}x^m y^n,\ H(x):=\frac {\psi_1(x)\chi_{_E}(x) +\psi_2(x) \chi_{_J}(x)}{S_b(x)}.
\end{split}
\ee
From \eqref{jumpG0_a}, \eqref{s_pm}, and \eqref{abeq} it follows that $H(z)$ is analytic in a neighborhood of any $z\in\mathring E\cup \mathring J$.

Let $R_+(x,y,\l)$ and $R_-(x,y,\l)$ be the analytic continuations in $\l$ of the kernel of the resolvent across the cut $[-1,1]$ from above and from below, respectively. It follows from Lemma~\ref{propRS} that $R_\pm(x,y,\l)$ are meromorphic functions of $\l$ for any $x,y$ that do not coincide with an endpoint of $J$ and $E$. The locations of the poles, of course, are independent of the choice of $x,y$. Let $\l_0$ be a pole of, say, $R_+(x,y,\l)$. Then, by symmetry (see Remark \ref{rem-R-pm}), $\l_0$ is also a pole of $R_-(x,y,\l)$. In what follows, with a slight abuse of notation, we denote by $R(x,y;\l)$ the kernel, which is the average of $R_+(x,y,\l)$ and $R_-(x,y,\l)$. The residue of $R(x,y;\l)$ at $\l_0$ is then:
\be\label{res-def}
\res{\l=\l_0} R(x,y;\l)=\frac12\left(\res{\l=\l_0} R_+(x,y;\l)+\res{\l=\l_0} R_-(x,y;\l)\right).
\ee
Applying \eqref{P-symm} to $R_+(x,y,\l)$ and $R_-(x,y,\l)$, combining the two residues, and using the symmetry of $R_\pm$ (see Remark~\ref{rem-R-pm}) we get that the residue of $R$ defined in \eqref{res-def} equals
\be
\label{P-symm-sum}
P_{\l_0}(x,y) =\Re \left(B(x,y) H(x)H(y)\right).
\ee
Clearly, $P_{\l_0}(x,y)$ in \eqref{P-symm-sum} is real-valued and satisfies
$P_{\l_0}(x,y)=P_{\l_0}(y,x)$.




\bt
Let $\l_0$ be an eigenvalue of  $\mathscr K$ imbedded in the continuous spectrum. Then the corresponding eigenspace  has a finite dimension bounded by twice the degree of the rational function $h(z)$ from Proposition \ref{prop-pole}. 
\et

{\bf Proof.} First, we show that the residue of the resolvent defined according to \eqref{res-def} defines the projector in $L^2(U)$ onto the corresponding eigenspace. A similar statement in the case of an isolated eigenvalue is well-known \cite{ReedSimon1}. Here our situation is a bit more complex, since all the eigenvalues are imbedded in the continuous spectrum. Nevertheless, the proof is fairly straightforward. We could not find a reference in any of the well-known texts on operator theory, so we decide to give it here for completeness. As is known, the projector onto the eigenspace of $\scr K$ corresponding to $\l_0$ can be computed as follows:
\bea\label{ResOfIdGeneral-v2}
\mathscr P_{\l_0}=\lim_{\delta\to0}\hat{E}_{(\l_0-\delta,\l_0+\delta)}=\frac{-1}{2\pi i\l_0}\lim_{\delta\to0}\lim_{\epsilon\to0^+}\int_{\l_0-\delta}^{\l_0+\delta}\left[
    \mathscr{R}(t+i\epsilon)-\mathscr{R}(t-i\epsilon)\right]dt,
\eea
where all the limits are in the sense of strong operator convergence. Pick any two functions $\phi_{1,2}\in C_0^\infty(\mathring E\cup \mathring J)$. 
Using that the kernels of $\mathscr{R}_\pm$ are analytic with respect to $x,y$ away from the endpoints of $J$ and $E$ and is a meromorphic function of $\l$, it is immediate that
\be\begin{split}\label{step-one}
 \lim_{\epsilon\to0^+}\int_{\l_0-\delta}^{\l_0+\delta}
    (\mathscr{R}(t\pm i\epsilon)\phi_1,\phi_2)dt
&= \lim_{\epsilon\to0^+}\int_{\l_0-\delta}^{\l_0+\delta}
    (\mathscr{R}_\pm(t\pm i\epsilon)\phi_1,\phi_2)dt\\
&=\int_{B_\pm(\delta)}\mathscr{R}_\pm(\l)d\l=\pm\pi i \res{\l=\l_0} (\mathscr{R}_\pm(\l)\phi_1,\phi_2)+O(\delta),
\end{split}
\ee
where $B_\pm(\delta)$ is the half-circle centered at $\l_0$ with radius $\delta$ in the upper and
lower halfplanes, respectively, oriented in the counter clockwise direction. Substituting \eqref{step-one} into \eqref{ResOfIdGeneral-v2} we obtain
\bea\label{projector}
(\mathscr P_{\l_0}\phi_1,\phi_2)=\frac{-1}{2\pi i\l_0}\pi i\left((\res{\l=\l_0} \mathscr{R}_+(\l)+ \res{\l=\l_0}\mathscr{R}_-(\l))\phi_1,\phi_2\right)=\frac{-1}{\l_0}(\res{\l=\l_0} \mathscr{R}(\l)\phi_1,\phi_2),
\eea
where the last equality follows from the definition \eqref{res-def}. Therefore the operators on the left and on the right in \eqref{projector} act the same way on $C_0^\infty(\mathring E\cup \mathring J)$. Comparing \eqref{projector} and \eqref{res-def} implies that the kernel of $\mathscr P_{\l_0}$ is the expression in \eqref{P-symm-sum}.


Since $P_{\l_0}(x,y)$ is self-adjoint, real-valued (cf. \eqref{P-symm-sum}) and degenerate (cf. \eqref{P-symm}), we can represent it in the form
\be\label{re-part-zero}\begin{split}
P_{\l_0}(x,y) =\sum_{m,n=1}^{N}a_{mn}f_m(x)f_n(y)=\sum_{n=1}^{N'} a_n g_n(x)g_n(y)
\end{split}
\ee 
for some real, symmetric matrix $(a_{mn})$ and real-valued functions $f_n(x)\in C^\infty(\mathring E\cup \mathring J)$. The latter are are analytic on $\mathring E\cup \mathring J$. Here $a_n$ are non-zero eigenvalues of $(a_{mn})$ (hence, $N'<N$), and $g_n$'s are obtained by a unitary transformation from $f_n$'s. Without loss of generality we can assume that the set of functions $\{f_n\}$ is linearly independent in $C^\infty(\mathring E\cup \mathring J)$. Then the set $\{g_n\}$ is linearly independent as well.
To prove that $\mathscr P_{\l_0}$ is of finite rank we just need to show that $g_n\in L^2(U)$ for all $n$. Clearly, we can find $\phi\in C_0^\infty(\mathring E\cup \mathring J)$, such that $\int_U g_n(x)\phi(x)dx\not=0$ if $n=1$ and equals zero for all other $n$. Since $\mathscr P_{\l_0}\phi\in L^2(U)$, it follows that $g_1\in L^2(U)$. Repeating the same argument for all $n$ implies the desired result.

Finally, from \eqref{bezout}, \eqref{P-symm}, \eqref{P-symm-sum}, and \eqref{re-part-zero} it follows that the dimension of the eigenspace of $\scr K$ corresponding to $\l_0$ does not exceed twice the degree of the rational function $h$, i.e. $2 \max (\deg S_b, \deg S_t)$.
\QED

\subsection{Proof of Theorem~\ref{theo-mainK^2}}
By changing variables
\be
\scr K^2=\int \l^2 dE_\l=\int_{\infty}^0 t dE_{-\sqrt t}+\int^{\infty}_0 t dE_{\sqrt t}
=\int^{\infty}_0 t d(E_{\sqrt t}-E_{-\sqrt t}),
\ee
it follows that the resolution of the identity associated with $\mathscr K^2$ is given by:
\be\label{res-id-k2}
V_\l=\begin{cases}
E_{[-\sqrt\l,\sqrt\l]},&\l>0,\\
0,&\l \leq 0.\end{cases}
\ee
Here $E_{[a,b]}:=\text{s-lim}_{\delta\to0^+}E_{(a-\delta,b]}$, and $E_\l$ is assumed to be strongly continuous from the right, see \cite{Kato}, Section X.1.1. The above definition ensures that $V_\l$ is strongly continuous from the right as well.
Assertions 1, 4, and 5 as well as the first half of assertion 2 follow directly from \eqref{res-id-k2} and the corresponding assertions of Theorem \ref{theo-main}.

To prove the second half of assertion 2 we show that $A A^\dagger$ and $A^\dagger A$ are unitarily equivalent. Indeed, let $A=V(A^\dagger A)^{1/2}$ be the polar decomposition of $A$. Here $V$ is a partial isometry with $\text{Ran} V=\overline{\text{Ran} A}$, which is uniquely defined by the condition $\text{Ker} V=\text{Ker} A$ (see Section VI.7 in \cite{Kato}). Using that $\text{Ker} A^\dagger=(\text{Ran} A)^\perp$ (see eq. (5.10) in Chapter III of \cite{Kato}) and that both $(-1)A$ and $A^\dagger$ are Hilbert transforms (i.e., densely defined with zero kernels), it follows that $\text{Ker} V=\text{Ker} A=0$ and $\text{Ran} V=\overline{\text{Ran} A}=L^2(E)$. Hence $V:\,L^2(J)\to L^2(E)$ is an isometry. Then $A A^\dagger=VA^\dagger A V^\dagger$, and the result follows. Therefore, in particular, the absolutely continuous spectra of $A A^\dagger$ and $A^\dagger A$ have the same multiplicity, and the latter equals to half of the multiplicity of $\Sp_{ac}(\mathscr K^2)$.

The first half of assertion 3 is proven similarly to assertion 3 of Theorem~\ref{theo-main}. Alternatively, this statement can be proven by showing that if $\l^2$ is an eigenvalue of $B$, then $\l$ (and $-\l$ if $\l\not=0$) is an eigenvalue of $\scr K$, and then invoking Theorem~\ref{theo-main}.
If there are no double points, Theorem~\ref{theo-main} implies that $\mathscr K^2$ and, therefore, $B$ is of trace class. If there are double points, then $\mathscr K$ and $\mathscr K^2$ have continuous spectrum and cannot be trace class. The last statement of assertion 3 follows from the standard operator theory.
\QED

\br\label{rem-asser23}
In the proof of Theorem \ref{theo-mainK^2} we showed that $\Sp_{ac}(A A^\dagger)=\Sp_{ac}(A^\dagger A)=[0,1]$ with the same multiplicity $n$.
It is instructive to prove this assertion directly by  following the arguments of Theorem \ref{ThmAtrace}. Consider first the simple case of $r=1$ in \eqref{JE-defs}, with both endpoints of $U=U_1$ belonging to $E=E_1$ (Lemma~\ref{Lemma1}, first case). Using  \eqref{KOL}, we obtain
\be\label{K0^2}
\mathscr K_0^2 = \mathscr K_{\text{ext}}^2 + \mathscr S^2+\mathscr K_{\text{ext}} \mathscr S +  \mathscr S \mathscr K_{\text{ext}}.
\ee 
According to the proof of Lemma~\ref{Lemma1}, $\mathscr S$ and, therefore,
 $\mathscr S^2+\mathscr K_{\text{ext}} \mathscr S +  \mathscr S \mathscr K_{\text{ext}}$, are trace class operators. Hence the absolutely continuous parts of $\mathscr K_0^2$ and $\mathscr K_{\text{ext}}^2$ are unitarily equivalent. Since: (a) $\mathscr K_0^2$ and $\mathscr K_{\text{ext}}^2$ are block diagonal relative to the decomposition $L^2(\R)=L^2(\wh E)\oplus L^2(J)$; (b) diagonal blocks of a block-diagonal trace class operator are also trace class, and; (c) each diagonal block of $\mathscr K_0^2 $ has absolutely continuous spectrum $[0,1]$ of multiplicity $n$ (see \cite{kbt19}), we conclude that the absolutely continuous spectrum of each block of $\mathscr K_{\text{ext}}^2$ is the interval $[0,1]$, and its multiplicity equals $n$. Restricting the blocks of $\mathscr K_{\text{ext}}^2$ to the blocks of $\mathscr K^2$ similarly to how this is done at the end of the proof of the first case in Lemma~\ref{Lemma1}, we obtain the desired assertion.
 
In a similar fashion, we use \eqref{KKS} and \eqref{main-K_split} to prove the assertion in all the remaining cases. The key observation is that all the cross terms not containing trace class operators
are zero when the right-hand sides of \eqref{KKS} and \eqref{main-K_split} are squared. Consider, for example, \eqref{KKS}. Now $\mathscr K_0^2 $ is a trace class perturbation of $ (\mathscr K_{\text{ext}} +  \mathscr K_\infty )^2= \mathscr K_{\text{ext}}^2 +  \mathscr K_\infty^2$, because, by construction, $\mathscr K_{\text{ext}} \mathscr K_\infty =   \mathscr K_\infty \mathscr K_{\text{ext}}=0$.  In \eqref{main-K_split} we get from \eqref{Kj-ker}
\be\label{main-K_split-sq}
\left(\bigoplus _{j=1}^{ r } \mathscr K_j\right)^2=\bigoplus _{j=1}^{ r } \mathscr K_j^2,
\ee
and the rest of the argument is analogous.
\er

\newcommand{\bto}{\beta_{\text{od}}}
\newcommand{\bte}{\beta_{\text{ev}}}
\newcommand{\bet}{\beta} 
\newcommand{\btt}{\beta} 
\def\a{\alpha}
\newcommand{\bff}{{\bf f}} 
\newcommand{\CB}{\mathcal B} 
\newcommand{\CM}{\mathcal{M}}
\newcommand{\la}{\lambda}  
\newcommand{\aux}{\mathring{\phi}} 
\newcommand{\Tin}{T_{\text{in}}}
\newcommand{\Tex}{T_{\text{ex}}}
\newcommand{\ex}{\text{ex}}
\newcommand{\inn}{\text{in}}

\appendix

\section{Spectrum of the operator $\mathscr K$ when $U=\R$}
In this section we extend an approach, which was originally developed
in  \cite{kbt19}, see also \cite{kbt20}. One is given a collection of $2n$ points $b_j\in\R$, $1\leq j\leq 2n$ (i.e., all $b_j$ are double points). We assume that they are arranged in ascending order: $b_j<b_{j+1}$, $1\leq j<2n$. Define
\be\label{intervals-def}\begin{split}
\Iin&:=[b_1,b_2]\cup[b_3,b_4]\cup\dots\cup [b_{2n-1},b_{2n}], \\ \Iex&:=[b_2,b_3]\cup[b_4,b_5]\cup\dots\cup [b_{2n},b_1],\ [b_{2n},b_1]:=(-\infty,b_1]\cup [b_{2n},\infty).
\end{split}
\ee
In terms of \eqref{JE-defs} this means that $r=1$ and $U_1=U=\R$. 
We have assumed that the point at infinity belongs to $\Iex$, but this does not affect the generality of the argument due to Lemma~\ref{lemmamobius}. Define
\be\label{oepols}
\begin{split}
\bto(z)&= \prod_{j=1}^{n}(z-b_{2j-1}),\quad 
\bte(z)= \prod_{j=1}^{n}(z-b_{2j}), \quad \bet(z)=\bte(z)/\bto(z),\\ 
\aux(z)&=\ln \bet(z),\ \phi(z)=\Re \aux(z),
\end{split}
\ee
where  we choose the standard branch of the logarithm.


The following facts are proven for $x\in\Iin$ in \cite{kbt19}, and the proofs for $x\in\Iex$ are analogous:
\begin{enumerate}
\item We have 
\be\label{arg-ln}
\Im \aux(x)=\pi,\ x\in\Iin,\quad \Im \aux(x)(x)=0,\ x\in\Iex;
\ee
\item The behavior of $\phi$ on the subintervals $(b_{2j-1},b_{2j})\subset\Iin$ and $(b_{2j},b_{2j+1})\subset\Iex$ satisfies
\be\begin{split}\label{phi-beh}
&\phi'(x)<0, x\in \Iin,\ \phi(x)\to +\infty,\ x\to b_{2j-1}^+,\ \phi(x)\to -\infty,\ x\to b_{2j}^-,\\
&\phi'(x)>0, x\in \Iex,\ \phi(x)\to -\infty,\ x\to b_{2j}^+,\ \phi(x)\to +\infty,\ x\to b_{2j+1}^-;
\end{split}
\ee
Therefore, $\phi(x)$ is monotonic and invertible on each subinterval, and the range of $\phi(x)$  on each subinterval is $\R$;
\item One has
\be\label{der-2}
\phi'(x)=\frac{Q(x)}{\bto(x)\bte(x)},\quad
Q(x):=\bte'(x)\bto(x)-\bte(x)\bto'(x),
\ee
and $Q(x)>0$ is bounded away from zero on $\R$;
\end{enumerate}
Suppose $s=\phi(z)/2$, $z\in(b_{2m},b_{2m+1})\subset\Iex$  and $t=\phi(x)/2$, $x\in(b_{2k-1},b_{2k})\subset\Iin$. Then
\be\label{cos_sinh_1}
\begin{split}
\cosh(s-t)&=\cosh\le( \frac{\phi(z)-\phi(x)}2\ri)= \cosh\le( \frac{\aux(z)-(\aux(x)-i\pi)}2\ri)
= i\sinh\le( \frac{\aux(z)-\aux(x)}2\ri).
\end{split}
\ee
Moreover,
\be\label{cos_sinh_2}
\begin{split}
2\sinh\le( \frac{\aux(z)-\aux(x)}2\ri)&= \sqrt{\frac{\bet(z)}{\bet(x)}}-\sqrt{\frac{\bet(x)}{\bet(z)}}=\frac{\bte(z)\bto(x)- \bte(x)\bto(z) }
{D}\\
&=\frac{(z-x)\sum_{i,j=1}^n B_{ij}z^{i-1}x^{j-1}}{D}
=\frac{(z-x)\sum_{j=1}^{2n} \rho_j P_j(z)P_j(x)}{D}, 
\end{split}
\ee
where $B:=B(\bte,\bto)=(B_{ij})$ is the B\'{e}zout matrix of the polynomials $\bte(z),\bto(z)$, and
\be\label{D-deff}
D:=\bto(z)\sqrt{\bet(z)}\bto(x)\sqrt{\bet(x)}.
\ee
Using \eqref{arg-ln} we find
\be\label{D-full}\begin{split}
D&=\text{sgn}(\bto(z))i\sqrt{\prod_{j=1}^{2n}|z-b_j|}\,\text{sgn}(\bto(x))\sqrt{
\prod_{j=1}^{2n}|x-b_j|}\\
&=i\text{sgn}(\bto(z))\text{sgn}(\bto(x))\sqrt{
\prod_{j=1}^{2n}|x-b_j||z-b_j|}.
\end{split}
\ee
Introduce two isometries
\be\label{tr-three}
\begin{split}
&\Tin:\ L^2(\Iin)\to L_n^2(\R),\ \Tex:\ L^2(\Iex)\to L_n^2(\R),\\
&\check\bff_\inn (t):=(\Tin f)(t):=\sqrt2\left(\left.\frac{{\text{sgn}(\bto(x))}f(x)}{\sqrt{|\phi'(x)|}}\right|_{x=\phi_1^{-1}(2t)},\dots,\left.\frac{{\text{sgn}(\bto(x))}f(x)}{\sqrt{|\phi'(x)|}}\right|_{x=\phi_{2n-1}^{-1}(2t)}\right),\\
&\check \bff_\ex(s):=(\Tex f)(s)
=\sqrt2\left(\left.\frac{{\text{sgn}(\bto(z))}f(z)}{\sqrt{\phi'(z)}}\right|_{z=\phi_{2}^{-1}(2s)},\dots,\left.\frac{{\text{sgn}(\bto(x))}f(z)}{\sqrt{\phi'(z)}}\right|_{z=\phi_{2n}^{-1}(2s)}\right),
\end{split}
\ee
where $L_n^2(\mathbb R)$ is the direct sum of $n$ copies of $L^2(\mathbb R)$, $L^2_n(\mathbb R) =\oplus_{j=1}^{n} L^2(\mathbb R)$. Here we set $\Vert \check\bff\Vert^2=\Vert\check f_1\Vert^2+\dots+\Vert\check f_n\Vert^2$, where $\check\bff=(\check f_1,\dots,\check f_n)\in L_n^2(\mathbb R)$ and $\Vert\check f_m\Vert$ is the conventional $L^2(I_m)$ norm. Also, in \eqref{tr-three}, $\phi_k^{-1}$ is the inverse of $\phi(x)$ on the $k$-th interval $(b_k,b_{k+1})$. By convention, the $2n$-th interval is $\R\setminus (b_1,b_{2k})$, i.e. it includes the point at infinity.

Changing variables in the definition of $A$ gives
\be\label{change-1}
\begin{split}
(\Tex A \Tin^{-1}\check \bff_\inn)_m(s)
&=\frac{\text{sgn}(\bto(z_m))}{\pi}\sqrt{\frac2{\phi'(z_m)}}\sum_{k=1}^n \int_{\mathbb R}\frac{\text{sgn}(\bto(x_k))\check f_k(t)}{\sqrt{|\phi'(x_k)|/2}\,(z_m-x_k)}dt\\
&=\frac{2\text{sgn}(\bto(z_m))}{\pi}\sum_{k=1}^n \int_{\mathbb R}\frac{\text{sgn}(\bto(x_k))\check f_k(t)}{\sqrt{|\phi'(x_k)|\phi'(z_m)}\,(z_m-x_k)}dt,\\
x_k:&=\phi_k^{-1}(2t),\,z_m:=\phi_m^{-1}(2s).
\end{split}
\ee
Combining \eqref{change-1},  \eqref{der-2}, \eqref{cos_sinh_1}, and \eqref{cos_sinh_2} we find
\be\label{xmz}
\frac{\text{sgn}(\bto(x_k))\text{sgn}(\bto(z_m))}{\sqrt{|\phi'(x_k)|\phi'(z_m)}\,(z_m-x_k)}=\frac{1}{2\cosh(s-t)}
\sum_{j=1}^n \frac{\rho_j P_j(x_k)P_j(z_m)}{\sqrt{Q(x_k)Q(z_m)}}.
\ee
Define two matrix functions
\be\label{one-matr}
\begin{split}
\CM_\inn:&=\{M_{jk}^{(\inn)}(t)\},\ M_{jk}(t):=P_j(x_k)\sqrt{\frac{\rho_j}
{Q(x_k)}},\ x_k:=\phi_k^{-1}(2t),\\
\CM_\ex:&=\{M_{jm}^{(\ex)}(s)\},\ M_{jm}(s):=P_j(z_m)\sqrt{\frac{\rho_j}
{Q(z_m)}},\ z_m:=\phi_m^{-1}(2s).
\end{split}
\ee
It is shown in \cite{kbt19} that $\{M_{jk}^{(\inn)}(t)\}$ is an orthogonal matrix for all $t\in\R$. The proof that $\{M_{jk}^{(\inn)}(s)\}$, $s\in\R$, is an orthogonal matrix is analogous. Substituting \eqref{xmz} and \eqref{one-matr} into \eqref{change-1} gives
\be\label{change-simple}
\begin{split}
(\Tex A \Tin^{-1}\check \bff_\inn)_m(s)
&=\sum_{j=1}^n  M_{jm}^{(\ex)}(s) \sum_{k=1}^n \int_{\mathbb R}\frac{M_{jk}^{(\inn)}(t) \check f_k(t)}{\pi\cosh(s-t)}dt.
\end{split}
\ee
In compact form, \eqref{change-simple} can be written as follows
\be\label{change-alt}
\Tex A \Tin^{-1}\check \bff_\inn = \CM_\ex^T K \CM_\inn\check \bff_\inn,
\ee
where $K$ is the operator of component-wise convolution with $(\pi\cosh(t))^{-1}$. 

Equation \eqref{change-simple} matches with the results in \cite{kbt19} in the case $n=1$ (see eq. (2.12) in \cite{kbt19}). Indeed, suppose $\Iin=[-b,b]$, and $\Iex=(-\infty,b]\cup [b,\infty)$. Then \eqref{der-2} and \eqref{cos_sinh_2} imply that $\rho_1=2b$ and $Q(x)\equiv 2b$, i.e. $\CM_\ex\equiv\CM_\inn\equiv1$ in \eqref{one-matr}. Observe also that there are two sign changes between \eqref{change-simple} and (2.12) in \cite{kbt19}. The first one arises because the operator $A$ in \eqref{Aoper} is negative of the Hilbert transform. The second sign change arises because $\Tex$ in \eqref{tr-three} is the negative of $\Tex$ in (2.11) of \cite{kbt19}. As a result, both in \eqref{change-simple} and in (2.12) of \cite{kbt19}, the corresponding operator becomes the convolution with  $(\pi\cosh(t))^{-1}$ after a change of variables.

Let $\CF:\, L_n^2(\mathbb R)\to L_n^2(\mathbb R)$ denote the map consisting of $n$ component-wise one-dimensional Fourier transforms (cf. \eqref{ft-inv}). Using \eqref{change-simple} and the integral 2.5.46.5 in \cite{pbm1},
we get 
\be\label{K-four}
K=\CF^{-1}\left(\frac1{\cosh(\pi\la/2)}\Id_n\right)\CF,
\ee
where $\la$ is the spectral (Fourier) variable, and $\Id_n$ is the $n\times n$ identity matrix. 
Therefore, \eqref{change-simple} gives
\be\label{H-final}
A f = (\CF\CM_\ex \Tex)^{-1} \left(\frac1{\cosh(\pi\la/2)}\Id_n\right) (\CF\CM_\inn \Tin)f.
\ee
Applying the adjoint to \eqref{H-final}, we get that $\mathscr K$ satisfies
\be\label{K-factor}\begin{split}
\mathscr K
=&\begin{bmatrix} 0 & A\\
A^\dagger & 0 \end{bmatrix}
=U^{-1}
\begin{bmatrix}
0 & \frac1{\cosh(\pi\la/2)}\Id_n\\
\frac1{\cosh(\pi\la/2)}\Id_n & 0 \end{bmatrix}U,\\
U:=&\begin{bmatrix}
\CF\CM_\ex \Tex & 0\\
0 & \CF\CM_\inn \Tin \end{bmatrix}: L_{2n}^2(\mathbb R)\to L_{2n}^2(\mathbb R),
\end{split}
\ee
where $U$ is an isometry. As is easily checked, the following self-adjoint isometry diagonalizes the middle operator on the right in \eqref{K-factor}
\be\label{Vdef}
V:=\frac1{\sqrt2}\begin{bmatrix} \Id_n & \Id_n\\
\Id_n & -\Id_n \end{bmatrix}: L_{2n}^2(\mathbb R)\to L_{2n}^2(\mathbb R),
\ee
therefore
\be\label{K-factor-v2}\begin{split}
\mathscr K
=&U^{-1}V^{-1}
\begin{bmatrix}
\frac1{\cosh(\pi\la/2)}\Id_n & 0\\
0 & -\frac1{\cosh(\pi\la/2)}\Id_n \end{bmatrix}VU.
\end{split}
\ee
The range of the function $(\cosh(\pi\la/2))^{-1}$ is $(0,1]$, and each value is taken twice. Hence we proved the following result.

\bt\label{exact-soln} Suppose $r=1$ in \eqref{JE-defs}, and $\Iin_1\cup \Iex_1=U_1=\R$, i.e $U$ consists of only one interval and coincides with all of $\R$. In this case the spectral interval of $\mathscr K$ is $[-1,1]$, the spectrum is absolutely continuous (i.e., there is no point spectrum), and its multiplicity equals to the number of double points (which is twice the number of subintervals in $\Iin_1$ or $\Iex_1$).
\et


\end{document}